\documentclass[10pt,psamsfonts]{amsart}

\usepackage{indentfirst}
\usepackage[dvips]{color}
\usepackage{graphics,amsthm}
\usepackage{color}
\usepackage{epsfig}
\usepackage{fancybox}
\usepackage{mathrsfs}
\usepackage{hyperref}
\usepackage{newlfont}
\usepackage{amsmath,amsfonts,amstext,amssymb,bezier,amsthm}
\usepackage{subfigure}
\usepackage[active]{srcltx}
\usepackage[all]{xy}
\usepackage[english]{babel}
\usepackage[latin1]{inputenc}
\usepackage{graphicx,psfrag}
\usepackage{fancyhdr}
\usepackage{niceframe}
\usepackage{xspace,setspace} %sum\'{a}rio
\usepackage{dsfont}
\usepackage{tikz}
\usetikzlibrary{matrix}
\usepackage{caption}%[font=small,labelfont=bf]{caption}

\newtheorem{theorem}{Theorem}[section]

\newtheorem{definition}{Definition}[section]
\newtheorem{proposition}{Proposition}[section]
\newtheorem{corollary}{Corollary}[section]
\newtheorem{lemma}{Lemma}[section]
\newtheorem{remark}{Remark}[section]

\newcommand*\xbar[1]{%
  \hbox{%
    \vbox{%
      \hrule height 0.5pt % The actual bar
      \kern0.5ex%         % Distance between bar and symbol
      \hbox{%
        \kern-0.1em%      % Shortening on the left side
        \ensuremath{#1}%
        \kern-0.1em%      % Shortening on the right side
      }%
    }%
  }%
}

\newcommand{\prooff}{\noindent \textbf{Proof: }}
\newcommand{\cqd}{{{\hfill $\square$}}}%\vspace{0.3cm}}

\newcommand{\R}{\mathbb{R}}

\newcommand{\I}{\textbf{I}}
\newcommand{\J}{\textbf{J}}
\newcommand{\K}{\textbf{K}}
\newcommand{\PP}{\textbf{P}}

\makeindex

\begin{document}
\renewcommand{\arraystretch}{1.5}

\author[Franzosa]
{Robert Franzosa$^1$}

%\author[de Rezende]{K. A. de Rezende$^2$}
%\thanks{$^2$Partially supported by CNPq under grant 302592/2010-5 and FAPESP under grant 2012/18780-0}

\author[Vieira]{Ewerton R. Vieira$^2$}
\thanks{$^2$Supported by FAPESP under grant 2010/19230-8.}

\address{$^1$ Departament of Mathematics and Statistics,
University of Maine, Orono, Maine, USA} \email{robert$\_$franzosa@umit.maine.edu}

\address{$^2$ Departamento de Matematica, Universidade
Estadual de Campinas, 13083--859, Campinas, SP,
Brazil} \email{ewertonrvieira@gmail.com}

%\address{$^3$ Departamento de Matematica, Universidade
%Estadual de Campinas, 13083--859, Campinas, SP,
%Brazil} \email{ewertonrvieira@gmail.com}

\title[Transition Matrix Theory]
{Transition Matrix Theory}

\subjclass[2010]{Primary 37B30; 37D15;  Secondary 70K70; 70K50; 55T05}

\keywords{Conley index, Connection Matrices, Transition Matrices, Morse-Smale System, Fast-Slow Systems}

\maketitle

%\author[ K. A. Rezende, R. Franzosa, E. R. Vieira]
%{Ketty A. de Rezende$^1$, Robert Franzosa$^2$ and Ewerton R. Vieira$^1$}

%\address{$^1$ Departamento de Matematica, Universidade
%Estadual de Campinas, Caixa Postal 6065, 13083--859, Campinas, SP,
%Brazil} \email{ketty@ime.unicamp.br, ewertonrvieira@gmail.com}

%\address{$^2$ Departament of Mathematics and Statistics,
%University of Maine, 0Orono, Maine, USA} \email{robert$\_$franzosa@umit.maine.edu}

%\subjclass[2010]{Primary 37B30; 37D15;  Secondary 70K70; 70K50; 55T05}

%\keywords{Conley index, Connection Matrices, Transition Matrices, Morse-Smale System; Sweeping Method; Spectral Sequence}

\maketitle

\begin{abstract}
%In this article we present the unification of the theory of algebraic, singular, topological and directional transition matrices by introducing generalized transition matrix which encompasses all of the previous four. Existence results are presented as well as bifurcation connections properties that all transition matrices share in common. By using those similarly features, important results were able to obtain by combining transition matrices, specially topological and singular transition matrix where both translate the interrelationship between singular dynamics and topological Conley index.

In this article we present a unification of the theory of algebraic, singular, topological and directional transition matrices by introducing the (generalized) transition matrix which encompasses each of the previous four. Some transition matrix existence results are presented as well as verification that each of the previous transition matrices are cases of the (generalized) transition matrix. Furthermore we address how applications of the previous transition matrices to the Conley Index theory carry over to the (generalized) transition matrix.
\begin{center}
Dedicated to the memory of James Francis Reineck
\end{center}
\end{abstract}
%$$
%\mathscr{ABCDEFGHIJLMNOPQRSTUVXZWY}
%$$
%$$
%\mathfrak{ABCDEFGHIJLMNOPQRSTUVXZWY}
%$$
%$$
%\mathcal{ABCDEFGHIJLMNOPQRSTUVXZWY}
%$$

\tableofcontents

\section{Introduction}

The Conley index theory has been a valuable topological technique for detecting global bifurcations in dynamical systems \cite{C}, \cite{CF}, \cite{FiM}, \cite{Fr1}, \cite{Fr2}, \cite{Fr3} and \cite{MM1}. This index is a standard tool in the analysis of invariant sets in dynamical systems, and its significance owes partly to the fact that it is invariant under local perturbation of a flow (the continuation property). Typically, one does not investigate a single invariant set in a dynamical system but rather works with decompositions of a larger invariant set into invariant subsets and connecting orbits between them. The Morse decomposition is the standard such decomposition in the Conley index theory. Within the index theory there are matrices of maps defined between the Conley indices of invariant sets in a Morse decomposition, and these matrices (the connection matrices) provide information about connections that exist between sets in the decomposition. The connection matrices also have local invariance properties under continuation. Nevertheless, under global continuation sets of connection matrices can undergo change which usually means that the dynamical system has undergone a global bifurcation.

An initial approach to identify bifurcations via the interplay between local invariance and global change in connection matrices was due to Reineck in \cite{R}(1988). By introducing an artificial slow flow on the parameter space in a continuous family of dynamical systems, he obtained a map between the Conley indices of Morse decomposition invariant sets at the initial parameter value and those at the final parameter value in a continuation. This map is known as a singular transition matrix, and it has the feature that a nonzero entry can identify a change of the connecting orbit structure of the Morse decomposition under the continuation.

Following Reineck's work on singular transition matrices, Mischaikov and McCord in \cite{MM1}(1992) used the continuation property of the Conley index, without introducing an artificial slow flow, to define matrices of maps between the Conley indices of Morse decomposition invariant sets at the initial and final parameter values in a continuation. Their maps, known as topological transition matrices, are naturally-defined maps on Conley indices that arise from the topological structure of invariant sets in the flow on the larger ``phase-cross-parameter" space. In order to define these maps on the indices of the Morse decomposition sets, it was necessary to assume that connection matrices were trivial at the end parameter values. Nonetheless, as with the singular transition matrices, they were able to demonstrate that non-trivial topological transition matrix entries identify potential bifurcations that exist in the overall continuation. Furthermore, in \cite{MM2} they established an equivalence between the singular transition matrices and the topological transition matrices in instances where both are defined. Later Franzosa, de Rezende and Vieira \cite {FdRV}(2014)  defined a new (general) topological transition matrix that extends the previous one, not requiring the assumption that the connection matrices are trivial at the end parameter values. In their case, the general topological transition matrix is defined to cover naturally-defined Conley-index maps rather than being defined directly by them.

%Following next, Mischaikov and McCord in \cite{MM1}(1992) use the continuation property of the Conley index without introducing an artificial slow flow, notwithstanding for some technical difficulties they had to assume the nonexistence of connection at the initial and final parameters in the parameter space. Even if the Conley index stays invariant under a continuation, its generators might be changed by the flow-defined Conley index isomorphism. Ensuing, the authors defined  topological transition matrix as map (between Conley index graded modules of the initial and final parameter) capable to track thoses changes of generators. Later Franzosa, de Rezende and Vieira \cite{FdRV}(2014) removed the assumption of no connection at the initial and final parameters by defined a new topological transition matrix that extends the previous one. In this new definition, the connection matrix may change the generators of the Conley index of the minimal Morse sets, which contrasts the previous case when no connections imply trivial connection matrix. Both cases provide global bifurcation when the entry is nonzero.

In \cite{FM}(1995), Franzosa and Mischaikov introduced the concept of an algebraic transition matrix. Given that connection matrices for a Morse decomposition are not unique, they raised the question of whether the nouniqueness could be understood via similarity transformations between connection matrices. Such transformations are algebraically defined, and - besides being associated with nonuniqueness of connection matrices for a particular Morse decomposition - can be exploited to track changes in connection matrices under flow continuation. They developed an existence result for algebraic transition matrices under the assumption of a ``stackable" underlying partial order, and they demonstrated how algebraic transition matrices could also be used to identify global bifurcations in a dynamical system under continuation. 

In \cite{KMO} and \cite{GKMOR},  the authors developed the directional transition matrix, a transformation that is similar in nature to both the singular transition matrix and the topological transition matrix. As in the singular transition matrix case, a slow flow is added to the parameter space, but it is more general than the specific flow used in defining the singular transition matrix. The advantage to the more-general approach is that it allows us to detect broader families of bifurcation orbits under continuation than those that are detected by the singular and topological transition matrices. The directional transition matrix is a transformation between indices of Morse decomposition sets at each end of the continuation, but not simply from those at one end of the continuation to those at the other (as in the other types of transition matrices). Instead, it maps the indices of those sets on either end that have an outward slow-flow direction to the indices of those sets with an inward slow-flow direction. As with the classical topological transition matrix, it is assumed that on each end of the continuation there are no connecting orbits between the Morse sets, so that natural flow-defined maps can be used to define the directional transition matrix. And, as in each of the above cases, the authors demonstrate how non-trivial directional transition entries identify bifurcations that occur under continuation.

While these four types of transition matrices are each defined differently and in different settings, they have in common that each is a Conley-index based algebraic transformation that tracks changes in index information under continuation and thereby identifies global bifurcations that could occur during the continuation. It is natural to expect that the theories could be unified in an overarching transition matrix theory, and that is the main purpose of this paper. The basic idea for this general transition matrix is that it covers natural flow-defined index isomorphisms that arise under a continuation. We discuss this aspect of the transition matrix further at the start of the next section, and then we provide our definition of the transition matrix. First, though,  we present briefly the necessary background material from Conley Index Theory, Morse decompositions, homology index braids, connection matrices, etc. (see \cite{C}, \cite{Fr1}, \cite{Fr2}, \cite{Fr3}, \cite{MMr} and \cite{S}).

%For instance, in \cite{FM}, is is proved that algebraic transition matrices can be used to generate all connection matrices at a particular parameter value.

%By undertaking the task of unifying the theory of transition matrix we actually present a framework which makes possible the defining of a generalized transition matrix which encompasses all transition matrices previously defined, such as those in \cite{R}, \cite{MM1} and \cite{FM}.

Throughout the paper $\PP$ represents a finite set with a partial order $<$. An \textit{interval} in $\PP$ is a set $\I \subseteq \PP$ which is such that if $p,q \in \I$ and $p < r < q$ then $r\in \I$. The set of intervals in $<$ is denoted by $\I(<)$. 

%An \textit{adjacent n-tuple of intervals} in $<$ is an n-tuple of mutually disjoint nonempty intervals, $(\I_1,...,\I_n)$, whose union is an interval in $<$ and which are such that if $p \in \I_i,$ $q \in \I_j$, and $p < q$, then $i < j$. 

An \textit{adjacent n-tuple of intervals} in $<$ is an ordered collection  $(\I_1,...,\I_n)$ of mutually disjoint nonempty intervals in $<$ satisfying: 
\begin{itemize}
\item $\bigcup_{i=1}^n \I_i\in \I(<);$
\item $\pi\in\I_j, \pi'\in\I_k,$ $j<k$ imply $\pi'\nless \pi.$
\end{itemize}

The collection of adjacent $n$-tuples of intervals in $<$ is denoted $\I_n(<)$. An adjacent 2-tuple of intervals is also called an \textit{adjacent pair} of intervals. If $<'$ is either an extension of $<$ or the restriction of $<$ to an interval in $<$, then $\I_n(<')\subseteq \I_n(<)$. If $(\I,\J)$ is an adjacent pair (2-tuple) of intervals, then $\I\cup\J$ is denoted $\I\J$. If $(\I_1,\ldots , \I_n)$ is called a decomposition of $\I$.

\begin{definition}
A graded module braid over $<$ is a collection $\mathcal{G}=\mathcal{G}(<)$ of graded modules and maps between the graded modules satisfying:
\begin{enumerate}
\item for each $\I\in\I(<)$, there is a graded modules $G(\I)$,
\item for each $(\I,\J)\in \I_2(<)$, there are maps: $i(\I,\I\J):G(\I)\rightarrow G(\I\J)$ of degree 0; $p(\I\J,\J):G(\I\J)\rightarrow G(\J)$ of degree 0; $\partial(\J,\I):G(\J)\rightarrow G(\I)$ of degree 1, that satisfy:

\begin{itemize}
\item $\cdots \rightarrow G(\I)\stackrel{i}{\rightarrow} G(\I\J)\stackrel{p}{\rightarrow} G(\J)\stackrel{\partial}{\rightarrow} G(\I)\stackrel{}{\rightarrow} \cdots$ is exact,
\item if $\I$ and $\J$ are noncomparable, then $p(\J\I,\I)i(\I,\I\J)=id_| G(\I)$,
\item if $(\I,\J,\K)\in\I_3(<)$, then the following braid diagram commutes. 

\end{itemize}
\end{enumerate}
\begin{center}
\includegraphics[width=8cm]{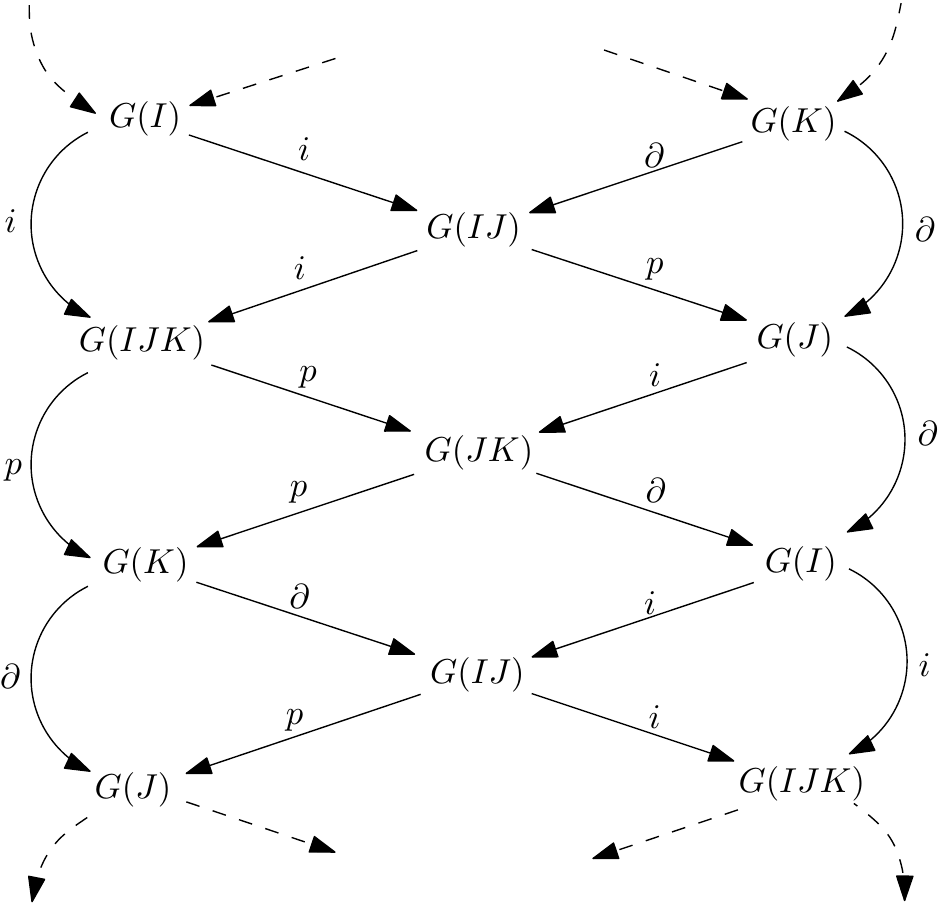}
\end{center}
\end{definition}

Assume that $\mathcal{G}$ and $\mathcal{G}'$ are graded module braids over $<$.

\begin{definition}
an $r$-degree map $\theta:\mathcal{G} \rightarrow \mathcal{G}'$ is a collection $\{\theta(\I)\}_{\I\in I(<)}$ of module homomorphisms $\theta(\I):G(\I)\rightarrow G'(\I)$ such that the following diagram commutes for each $(\I,\J)\in\I_2(<):$
$$
\xymatrixcolsep{3pc}\xymatrix{
\cdots \ar[r] & G_k (\textbf{I}) \ar[r]^{i} \ar[d]^{\theta(\I)} & G_k (\textbf{IJ}) \ar[r]^{p} \ar[d]^{\theta(\I\J)} & G_k (\textbf{J}) \ar[r]^{\partial_\lambda(\textbf{J},\I)} \ar[d]^{\theta(\J)}& G_{k-1}(\textbf{I})\ar[r] \ar[d]^{\theta(\I)}& \cdots\\
\cdots \ar[r] & G'_{k-r} (\I) \ar[r]^{i} & G'_{k-r} (\I\J) \ar[r]^{p} & G'_{k-r} (\J) \ar[r]^{\partial_\mu(\J,\I)} & G'_{k-r-1} (\I)\ar[r]& \cdots
}
$$

If, futhermore, $\theta(\I)$ is an isomorphism for each $\I\in \I(<)$, then $\theta$ is called an $r$-degree isomorphism and $\mathcal{G}$ and $\mathcal{G}$ are said to be $r$ isomorphic. 
\end{definition}

\begin{definition}
A chain complex braid over $<$ is a collection $C=C(<)$ of chain complexes and chain maps satisfying:
\begin{enumerate}
\item for each $\I\in\I(<)$, there is a chain complex $C(\I)$,
\item for each $(\I,\J)\in \I_2(<)$, there are $0$ degree maps $i(\I,\I\J):C(\I)\rightarrow C(\I\J)$ and $p(\I\J,\J):C(\I\J)\rightarrow C(\J)$ which satisfy:

\begin{itemize}
\item $ C(\I)\stackrel{i}{\rightarrow} C(\I\J)\stackrel{p}{\rightarrow} C(\J)$ is weakly exact,
\item if $\I$ and $\J$ are noncomparable, then $p(\J\I,\I)i(\I,\I\J)=id_| C(\I)$,
\item if $(\I,\J,\K)\in\I_3(<)$, then the following braid diagram commutes. 

\end{itemize}
\end{enumerate}
\begin{center}
\includegraphics[width=7cm]{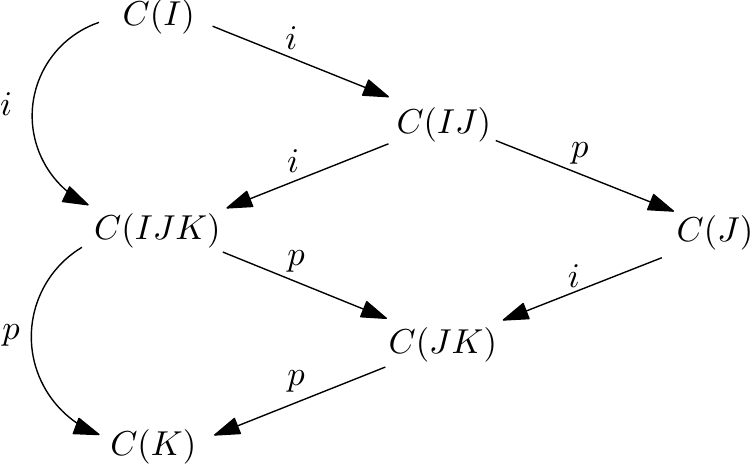}
\end{center}
\end{definition}

Now assume that $C$ and $C'$ are chain complex braids over $<$.

\begin{definition}
an $r$-degree chain map $T:C\rightarrow C'$ is a collection of maps $T(\I):C(\I)\rightarrow C'(\I), \I\in\I(<)$, such that for each $(\I,\J)\in\I_2(<)$ the following diagram commutes:
$$
\xymatrixcolsep{3pc}\xymatrix{
 C_k (\textbf{I}) \ar[r]^{i} \ar[d]^{T(\I)} & C_k (\textbf{IJ}) \ar[r]^{p} \ar[d]^{T(\I\J)} & C_k (\textbf{J})  \ar[d]^{T(\J)}\\
C'_{k-r} (\I) \ar[r]^{i} & C'_{k-r} (\I\J) \ar[r]^{p} & C'_{k-r} (\J) 
}
$$\vspace{0.001cm}
\end{definition}

To simplify and to agree with previous definitions, we denote: a $0$ degree map $\theta:\mathcal{G} \rightarrow \mathcal{G}'$ by map $\theta$; $0$ degree isomorphism $\theta:\mathcal{G} \rightarrow \mathcal{G}'$ by isomorphism $\theta$; and $0$ degree chain map $T:C\rightarrow C'$ by chain map $T$.

Let $\varphi$ be a continuous flow on a locally compact Hausdorff space and let $S$ be a compact invariant set under $\varphi$. A \textit{Morse decomposition} of $S$ is a collection of mutually disjoint compact invariant subsets of $S$,
$$
\mathcal{M}(S)=\{M(p)\ |\ \pi\in\PP\}
$$
indexed by a finite set $\PP$, where each set $M(p)$ is called a \textit{Morse set}. A partial order $<$ on $\PP$ is called an \textit{admissible ordering} if for $x\in S\backslash \bigcup_{\pi\in\PP}M(p)$ there exists $p<q$ such that $\alpha(x)\subseteq M(p)$ and $\omega(x)\subseteq M(q).$ The flow defines an
admissible ordering of $M$, called the \textit{flow ordering} of $M$, denoted $<_F$, and such
that $M(\pi)<_F M(\pi') $ if and only if there exists a sequence of distinct elements
of $P: \pi = \pi_0,\ldots,\pi_n=\pi'$, where $C(M(\pi_j),M(\pi_{j-1}))$, the set of connecting orbit between $M(\pi_j)$ and $M(\pi_{j-1})$, is nonempty for each $j = 1, \ldots ,n$. Note that every admissible ordering of $M$ is an extension of $<_F$.

In the Conley theory one begins with the Conley index for {isolated invariant sets}, i. e.,  $S\subseteq X$ is an  \textit{isolated invariant set} if there exists a compact set $N\subseteq X$ such that $S \subseteq int N$ and 
$$
S=Inv(N,\varphi)=\{x\in N|\ \mathcal{O}(x) \subseteq N\}.
$$
The \textit{homological Conley index of} $S$, $CH_\ast(S)$ is the homology of the pointed space $(N\backslash L)$, where $(N,L)$ is an index pair for $S$. Setting $$M(\I)=\bigcup_{\pi\in\I}M(\pi)\cup\bigcup_{\pi,\pi'\in \I}C(M(\pi'),M(\pi)),$$ the Conley index of $M(\I)$, $CH_\ast(M(\I))$, in short $H_\ast(\I)$, is well defined, since $M(\I)$ is an isolated invariant set for all $\I\in\I(<)$.

Given $\mathcal{M}(S)$, a Morse decomposition of $S$, the existence of an admissible ordering on $\mathcal{M}(S)$ implies that any recurrent dynamics in $S$ must be contained within the Morse sets, thus the dynamics off the Morse sets must be gradient-like. For this reason, Conley index theory refers to the dynamics within a Morse set as local dynamics and off the Morse sets as global dynamics. We briefly introduce the connection matrix theory, which addresses this latter aspect.

%\begin{definition}\label{connection_matrix}
%Given $\mathcal{G}$, a graded module braid over $<$, and $\mathcal{C} = \{C(\pi)\}_{\pi\in\PP}$,
%a collection of graded modules, let $\Delta: \bigoplus_{\pi\in\PP} C(\pi) \rightarrow \bigoplus_{\pi\in\PP} C(\pi)$ be a $<$-upper
%triangular boundary map. If $ \mathcal{H}\Delta$, the graded module braid generated by $\Delta$, is isomorphic to $\mathcal{G}$, and $C(p)$ is isomorphic to $\mathcal{G}(p)$ then $\Delta$
%is called a connection matrix of $\mathcal{G}$.
%\end{definition}

\begin{definition}\label{connection_matrix}
Given $\mathcal{G}$, a graded module braid over $<$, and $C=\{C(p)\}_{p\in\PP}$, a collection of graded modules, let $\Delta: \bigoplus_{p\in\PP} C(p)\rightarrow \bigoplus_{p\in\PP}C(p)$ be a $<$-upper triangular boundary map. Then:
\begin{enumerate}
\item If $\mathcal{H},$ the graded module braid generated by $\Delta$, is isomorphic to $\mathcal{G}$, then $\Delta$ is called a $C-$\textbf{\textit{connection matrix}} of $\mathcal{G}$;
\item If, furthermore, $C(p)$ is isomorphic to $G(p)$ for each $p\in \PP$, then $\Delta$ is called a \textbf{\textit{connection matrix}} of $\mathcal{G}$.
\end{enumerate}
\end{definition}

To simplify notation, for $\I \in \I(<)$ we denote $\bigoplus_{\pi\in\I}C(\pi)$ by $C(\I)$, and the
corresponding homology module in $\mathcal{H}\Delta$ by $H(\I)$. In particular, the homology index braid of an admissible ordering of a Morse decomposition $\mathcal{G}=\{H_\ast(\I)\}_{\I\in\I(<)}$ is an example of a graded module braid. In this setting  a $<$-upper
triangular boundary map
$$\Delta:\displaystyle \bigoplus_{\pi\in \PP}CH_\ast (M(\pi)) \rightarrow \displaystyle \bigoplus_{ \pi\in  \PP}CH_{\ast-1} (M( \pi))$$
satisfying Definition \ref{connection_matrix} for $\mathcal{C}\Delta=\{CH_\ast(M(\pi))\}_{\pi\in\PP}$ is called the \textit{connection matrix for  a Morse decomposition}. Moreover, let $\mathcal{CM}(<)$ denote the set of all connection matrices for a given ($<$-ordered) Morse decomposition $\mathcal{M}(S)$.

{One of the key features in Conley theory is its invariance under continuation. Since the connection matrices for Morse decompositions
are algebraically derived from the homology Conley index braid, this seems to indicate that connecting orbits that persist over open sets in parameter space are identified by connection matrices. We now define Conley index continuation.}

Let $\Gamma$ be a Hausdorff topological space, $\Lambda$ a compact, locally contractible, connected metric space and $X$ a locally compact metric space. Assume that $X\times\Lambda\subseteq \Gamma$ is a local flow and $Z$ is a locally compact space. Let $\Pi_X: X\times \Lambda\rightarrow X$ and $\Pi_\Lambda:X\times \Lambda\rightarrow \Lambda$ be the canonical projection maps. See \cite{S} and \cite{Fr3}.

\begin{definition}
A parametrization of a local flow $X\subseteq \Gamma $ is a homeomorphism $\phi :Z\times \Lambda \rightarrow X$ such that for each $\lambda \in  \Lambda$, $\phi (Z\times \{\lambda\})$ is a local flow.
\end{definition}

Let $\phi:Z\times \Lambda \rightarrow X$ be a parametrization of a local flow $X$. Denote the restriction $\phi |_{(Z\times \{\lambda\})}$ by $\phi_\lambda$ and its image by $X_\lambda$.

\begin{lemma}\emph{[Salamon]}
For any compact set $N\subseteq X$ the set $\Lambda(N)=\{\lambda\in\Lambda \ |\ N\times \lambda$ is an isolating neighborhood in $X\times \lambda\}$ is open in $\Lambda$.
\end{lemma}

\begin{definition}
The space of isolated invariant sets is
$$
\mathscr{S}= \mathscr{S}(\phi)=\{S\times \lambda\ |\ \lambda\in \Lambda\ \text{and}\ S\times \lambda\ \text{is an isolated invariant compact set in}\ X\times \lambda\}.
$$
\end{definition}

For all compact sets $N\subseteq X$ define the maps $\varrho_N:\Lambda (N)\rightarrow \mathscr{S}$ and $\varrho_N(\lambda)=Inv(N\times \lambda)$. Then consider the topology on the space $\mathscr{S}$ generated by the sets $\{ \varrho_N(U)\ |\ N\subseteq X $ {compact}, $U\subseteq \Lambda(N)$ {open} $\}$.

A map $\gamma:\Lambda\rightarrow \mathscr{S}$ is called a \textit{section} of the space of isolated invariant sets if $\Pi_\Lambda\circ\gamma=id|_\Lambda$.

{We are interested in the situation where the homology index braids of admissible orderings of Morse decompositions at parameters $\lambda$ and $\mu$ are isomorphic. That is, it is not enough that a Morse decomposition continues over $\Lambda$, it must also continue with a partial order, more specifically:}

\begin{figure}[!h]
  \centering
  \includegraphics[width=10cm]{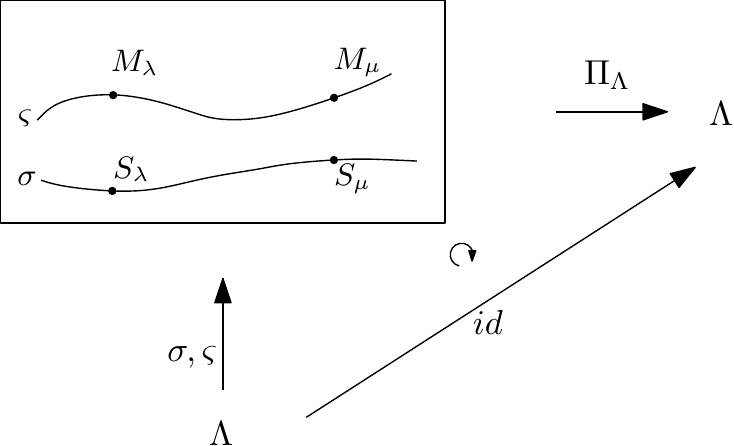}\\
  \caption{Sections from Definition \ref{def_continuation}}\label{secao}
\end{figure}

\begin{definition}\label{def_continuation}
Let $\mathcal {M}(S)=\{M(\pi)\ |\ \pi\in(\PP,<)\}$ be an ordered Morse decomposition of the isolated invariant set $S\subseteq X\times \Lambda$. Let $M_\lambda=\{M_\mu(\pi)\}_{\pi\in (\PP,<_\lambda)}$, $M_\mu=\{M_\mu(\pi)\}_{\pi\in (\PP,<_\mu)}$, $S_\lambda$ and $S_\mu$ be the sets obtained by intersection of $\mathcal{M}(S)$ and $S$ by the fibers $X\times\{\lambda\}$ and $X\times\{\mu\}$, respectively, where $<_\nu$ is the order restricted to the order $<$ in the parameter $\nu\in\Lambda$.

\begin{itemize}
  \item We say that $\mathcal{M}(S)$ with its order $<$ continues over $\Lambda$ if there exist sections $\sigma$ and $\varsigma_\pi : \Lambda \rightarrow \mathscr{S}$ such that $\{ \varsigma_\pi(\nu)|\ \pi \in (\PP,<_\nu)\}$ is a Morse decomposition for $\sigma(\nu),\ \forall \nu\in\Lambda$.
  \item  If, furthermore, there exist a path $\omega:[0,1]\rightarrow \Lambda$ from $\lambda$ to $\mu$; $\sigma(\lambda)=S_{\lambda}$; $\sigma(\mu)=S_{\mu}$; $\varsigma_\pi(\lambda)=M_{\lambda}(\pi)$; $\varsigma_\pi(\mu)=M_{\mu}(\pi)$; and if $\mathcal{M}(S)$ continues at least over $\omega([0,1])$, then we say that the admissible orderings $<_\lambda$ and $<_\mu$ are related by continuation or continue from one to the other. See Figure \ref{secao}.
\end{itemize}
\end{definition}

The following Lemma \ref{lemma_juncao} is helpful to understand the previous Definition \ref{def_continuation}.

\begin{lemma}\label{lemma_juncao}\emph{[McCord, Mischaikov, Salamon]}
\begin{itemize}
  \item Let $\gamma:\Lambda\rightarrow\mathscr{S}$ be a section, then $\gamma$ is continuous if and only if
  $$
  S=\displaystyle\bigcup_{\lambda\in\Lambda}\gamma(\lambda)
  $$
  is an isolated invariant set in $X\times \Lambda$.

  \item Let
  $$
  S=\displaystyle\bigcup_{\lambda\in\Lambda}\sigma(\lambda),\ \ \ M(\pi)=\displaystyle\bigcup_{\lambda\in\Lambda}\varsigma_\pi(\lambda)\ \text{for any}\ \pi\in \PP.
  $$
  Then, $S$ is an isolated invariant set in $X\times \Lambda$ under $\phi$ and $\mathcal{M}(S)=\{M(\pi)\ |\ \pi\in (\PP,<)\}$ is its Morse decomposition if, and only if, $\mathcal{M}(S)$ with its order continues.
\end{itemize}
\end{lemma}

Suppose that $S_0$ and $S_1$ are invariant sets related by continuation in $X_{\lambda_0}$ and $X_{\lambda_1}$. Hence, there exists a map $\omega:[0,1]\rightarrow \Lambda$ such that $\omega(0)=\lambda_0$ and $\omega(1)=\lambda_1$ and an isolated invariant set $S$ over $\omega(I)$ such that $S_{\lambda_i}=S_i$. The inclusion $f_i:X_{\lambda_i}\rightarrow X\times \omega(I)$ induces an isomorphism $CH_\ast(S_i)\stackrel{f_{i\ast}}{\longrightarrow} CH_\ast(S)$, where $CH_\ast(S_i)$ and $CH_\ast(S)$ indicates the Conley homology indices of $S_i$ in $X_{\lambda_i}$ and of $S$ in $X\times \omega(I)$, respectively. Thus, there is an isomorphism, called \textit{\textbf{Conley flow-define isomorphism}}
$$
F_\omega:CH_\ast(S_0) \stackrel{f^{-1}_{1\ast}\circ f_{0\ast}}{\longrightarrow} CH_\ast(S_1),
$$
that depends on the endpoint-preserving homotopy class $\omega$. If $\pi_1(\Lambda)=0$ then $F_\omega$ is independent of the path $\omega$ and one writes $F_{\lambda_1,\lambda_2}$ instead of $F_\omega$. The Conley flow-defined isomorphism is well-behaved with respect to composition of paths: $F_{\lambda,\lambda}=id$, $F_{\mu,\nu}\circ F_{\lambda,\mu}=F_{\lambda,\nu}$ and $F_{\lambda,\mu}=F^{-1}_{\mu,\lambda}$. For more details see \cite{MM2} and \cite{S}. {To simplify notation, we denote $CH_\ast(M_\nu(\I))=H_{\ast,\nu}(\I)$ or just $CH(M_\nu(\I))=H_{\nu}(\I)$, where $\I\in\I(<_\nu)$ and $\nu\in\{\lambda,\mu\}$}.

\section{Transition Matrix}

Given a Morse decomposition of an isolated invariant set $S$ in a flow, an associated homology index braid contains a significant amount of algebraic information about the overall structure of the Morse sets and connecting orbits in $S$. Generally speaking, the homology index braid is a relatively unmanagable structure to examine and analyze in order to draw conclusions about the structure of $S$. Overcoming this difficulty is part of the motiviation for the connection matrix. The connection matrix is a matrix of maps between the indices of the Morse sets that produces a graded module braid that is isomorphic to the homology index braid. In this way we can think of the connection matrix as ``covering" the homology index braid. Consequently, it is reasonable to expect that the connection matrix itself contains information from which the invariant set structure can be understood.

The transition matrices discussed in the introduction can be viewed similarly, and this perspective motivates our general transition matrix defined below. Given an invariant set and ordered Morse decomposition that continues over a parameter interval, there is a family of continuation-defined isomorphisms between indices of the invariant subsets that continue. If bifurcations occur during the continuation, then this will be reflected in a change in the flow-defined mappings in the homology index braids on either end of the continuation. The idea behind the general transition matrix is that it is a matrix of mappings between indices of Morse sets from either end of the continuation that ``covers'' the continuation-defined isomorphisms. By covering the continuation-defined isomorphisms, the transition matrix then assists in detecting change that occurs in the homology index braids, and thus it is reasonable to expect that the transition matrix contains information revealing change in the invariant set structure under continuation.

In this section we introduce our general definition of the transition matrix, and we address some straightforward properties of it. It is important to have the viewpoint that the transition matrix theory that we present here is an algebraic theory. In fact, the same could be said for the various previous versions of transition matrices. The primary motivation for introducing and studying transition matrices has been in applications to the Conley index theory, but there are important questions that need to be addressed regarding algebraic aspects of transition matrices. In section 3 we show how previous transition matrix existence results can be adapted to provide specific existence results for the transition matrix introduced here. In sections 4-7 we examine applications of the transition matrix in the Conley index theory, demonstrating how previous results using specific versions of the transition matrix carry over to our general setting.

%\textcolor{blue}{In this section, we explorer the algebraic properties about the transition matrix. Note that, these features is independent to work related to Conley index theory. So, it is possible to apply these algebraic tool in other context. However, so far, we only know application to Conley index theory, since this arose from a necessity to define a algebraic tool to analyze the homology index braid for Conley index associate to a partial ordered Morse decomposition. Application for Conley index theory is done in the next sections 4, 5, 6 and 7.}

\begin{definition}\label{def_cover}
Given chain complex braids $\mathcal{C}$ and $\mathcal{C}'$ and graded module braids $\mathcal{G}$ and $\mathcal{G}'$, an $r$-degree chain map $\mathcal{T}:\mathcal{C}\rightarrow \mathcal{C}'$ is said to \textbf{cover} an $r$-degree isomorphism $\theta$ (relative to $\Phi$ and $\Phi'$) if for all $\I\in \mathcal{I}(<)$, we have that the following diagram commutes
$$
\xymatrixcolsep{3pc} \xymatrix{
\mathcal{HC}(\I) \ar[r]^{\mathcal{T}_\ast(\I)} \ar[d]_{\Phi(\I)} & \mathcal{HC}'(\I) \ar[d]^{\Phi'(\I)}\\
\mathcal{G}(\I) \ar[r]^{\theta(\I)}                   & \mathcal{G}'(\I)
}
$$
where $\mathcal{T}_\ast(\I)$ is the homology map induced from the chain map $\mathcal{T}(\I)$, and $\Phi:\mathcal{HC}\rightarrow \mathcal{G}$ and $\Phi':\mathcal{HC}'\rightarrow \mathcal{G}'$ are graded module braid homology isomorphisms.
\end{definition}

\begin{definition}\label{TM}
If{, in {Definition} \emph{\ref{def_cover}},} $\mathcal{C}$ and $\mathcal{C}'$ arise from connection matrices $\Delta:\bigoplus_\PP C(p)\rightarrow \bigoplus_\PP C(p),\ \Delta':\bigoplus_\PP C'(p)\rightarrow \bigoplus_\PP C'(p)$, respectively, and $\mathcal{T}$ arises from a matrix $T:\bigoplus_\PP C(p)\rightarrow \bigoplus_\PP C'(p)$ then $T$ is called an $r$-degree (\textit{generalized}) \textbf{transition matrix} for $\Delta$ and $\Delta'$.
\end{definition}

Beyond this point in the paper, ``transition matrix'' refers to the transition matrix as defined here unless we specifically refer to a particular previous type such as ``algebraic transition matrix'', ``topological transition matrix'', etc. {Also we denote a $0$ degree transition matrix $T$ by transition matrix $T$.} 

\begin{theorem}\label{gen_prop}
Let $\Delta:\bigoplus_\PP C(p)\rightarrow \bigoplus_\PP C(p)$ and $\Delta':\bigoplus_\PP C'(p)\rightarrow \bigoplus_\PP C'(p)$ be connection matrices of $\mathcal{G}$ and $\mathcal{G}'$, respectively. Let $T:\bigoplus_\PP C(p)\rightarrow \bigoplus_\PP C'(p)$ be a transition matrix for $\Delta$ and $\Delta'$. Assume that $\Phi:\mathcal{HC}\rightarrow \mathcal{G}$ and $\Phi':\mathcal{HC}'\rightarrow \mathcal{G}'$ are graded module braid isomorphisms, and that $T$ covers an isomorphism $\theta$ (relative to $\Phi$ and $\Phi'$). Then the transition matrix $T$ satisfies the following properties:
\begin{description}
  \item[(i)] $ T \circ \Delta=\Delta'\circ  T$;
  \item[(ii)] $T(\{p\})=id$ and $T$ is upper triangular with respect to $<$;
  \item[(iii)] $T$ is an isomorphism;
  \item[(iv)] $T^{-1}$ covers $\theta^{-1}$ (relative to $\Phi'$ and $\Phi$). Moreover, suppose that $\mathcal{T}_1$ and $\mathcal{T}_2$ cover $\theta_1:\mathcal{G}\rightarrow \mathcal{G}'$ (relative to $\Phi$ and $\Phi'$) and $\theta_2:\mathcal{G}'\rightarrow \mathcal{G}''$ (relative to $\Phi'$ and $\Phi''$), and arise from $T_1$ and $T_2$, respectively. Then $\mathcal{T}_2\circ \mathcal{T}_1$ covers $\theta_2 \circ \theta_1$ (relative to $\Phi$ and $\Phi''$) and arises from $T_2\circ T_1$.
\end{description}
\end{theorem}

\prooff

(\textbf{i}) Since $\Delta(\I)$ and $\Delta'(\I)$ are boundary maps and $T(\I)$ is a chain map, we have that $T(\I)\circ \Delta(\I) = \Delta'(\I) \circ T(\I)$ for all $\I$.

(\textbf{ii}) Change $F$ for $\theta$ in the proof of Theorem 5 in \cite{FdRV}.

(\textbf{iii})By item (ii) we have that $T$ is the identity on the diagonal and upper triangular, thus $T$ is an isomorphism.

(\textbf{iv}) It follows directly by the composition of the maps. \cqd

\section{Algebraic Existence Results}\label{Existence}

%As indicated in the previous section, we present existence of transition matrix when the partial order satisfy some conditions. It is unknown if there exists a transition matrix which covers an isomorphism in the general case for partial order \textcolor{blue}{(not necessary related to Conley index theory)}. However, next sections we show that if the isomorphism is actually the Conley flow-defined isomorphism then there exists transition matrix.

{As indicated in the previous section, the transition matrix theory is an algebraic theory, and there are important questions to be addressed regarding the theory solely from the algebraic viewpoint. The most basic is the existence question: Given  connection matrices $\Delta$ and $\Delta'$, does there exist a transition matrix for them? The answer is unknown for this most general situation, however, the previous existence results for the specific types of transition matrices provide us with positive existence results under appropriate conditions. In this section, we present a few such results. %arising from the algebraic transition matrix, singular transition matrix, topological transition matrix, and general topological transition matrix theories.}

To begin, we consider the situation that is the algebraic basis behind the existence of topological transition matrices in \cite{MM1}. In this case we have trivial connection matrices $\Delta:\bigoplus_\PP C(p)\rightarrow \bigoplus_\PP C(p)$ and $\Delta':\bigoplus_\PP C'(p)\rightarrow \bigoplus_\PP C'(p)$. So, for each $\I\in \mathcal{I}(<)$ we have $\mathcal{HC}(\I) = \bigoplus_\I C(p)$ and $\mathcal{HC}'(\I) = \bigoplus_\I C'(p)$. Suppose that the associated graded module briad isomorphisms are $\Phi:\mathcal{HC}\rightarrow \mathcal{G}$ and $\Phi':\mathcal{HC}'\rightarrow \mathcal{G}'$. Then for each $\I\in \mathcal{I}(<)$, the corresponding isomorphisms are mappings $\Phi(\I):\bigoplus_\I C(p)\rightarrow \mathcal{G}(\I)$ and $\Phi'(\I):\bigoplus_\I C'(p)\rightarrow \mathcal{G}'(\I)$. We have the following{, straighforward,} transition matrix existence result for this case. The proof is straightforward.

\begin{theorem}
Let $\Delta$ and $\Delta'$ be connection matrices of $\mathcal{G}$ and $\mathcal{G}'$, respectively. Assume that $\Delta$ and $\Delta'$ are trivial and that $\Phi:\mathcal{HC}\rightarrow \mathcal{G}$ and $\Phi':\mathcal{HC}'\rightarrow \mathcal{G}'$ are graded module braid isomorphisms. If $\Theta:\mathcal{G} \rightarrow \mathcal{G'}$ is an isomorphism between graded module braids, then $T=\Phi^{-1}(\PP) \circ \Theta(\PP) \circ \Phi(\PP) :\bigoplus_\PP C(p)\rightarrow \bigoplus_\PP C'(p)$
is a transition matrix for $\Delta$ and $\Delta'$.
\end{theorem}

Next we address an existence result that is based on the algebraic transition matrix existence result in \cite{FM}. In this case the result holds where the underlying partial order is of a particular type. Let $(\I_1,...,\I_n)$ be an adjacent n-tuple of intervals whose union is $\I$, then $(\I_1,...,\I_n)$ is called a \textit{decomposition} of $\I$, if $\I,\J \subseteq \PP$ are disjoint then we say that $\I$ and $\J$ are \textit{noncomparable} if  neither $p<q$ nor $p<q$ for every $p\in \I,$ $q\in\J$ and we say that $\J$ is \textit{totally greater} than $\I$ if $p < q$ for every $p\in\I$ and $q\in \J$.
A partial order $<$ on $\PP$ is called \textit{stackable} if there is a decomposition of $\PP$, $(\I_1,...,\I_n)$, such that $<$ restricted to each $\I_i$ is trivial and such that if $i < j$ then $\I_j$ is totally greater than $\I_i$. Note that a trivial order is stackable and a total order is stackable.

\begin{definition}
Given chain complex braids $\mathcal{C}$ and $\mathcal{C}'$ and graded module braids $\mathcal{G}$ and $\mathcal{G}'$, a chain map $\mathcal{T}:\mathcal{C}\rightarrow \mathcal{C}'$ is said to \textbf{weakly cover} an isomorphism $\theta$ (relative to $\Phi$ and $\Phi'$) if there exists $\I\in \mathcal{I}(<)$ such that the following diagram commutes
$$
\xymatrixcolsep{3pc} \xymatrix{
\mathcal{HC}(\I) \ar[r]^{\mathcal{T}_\ast(\I)} \ar[d]_{\Phi(\I)} & \mathcal{HC}'(\I) \ar[d]^{\Phi'(\I)}\\
\mathcal{G}(\I) \ar[r]^{\theta(\I)}                   & \mathcal{G}'(\I)
}
$$
where $\mathcal{T}_\ast(\I)$ is the homology map induced by the chain map $\mathcal{T}(\I)$, $\Phi:\mathcal{HC}\rightarrow \mathcal{G}$ and $\Phi':\mathcal{HC}'\rightarrow \mathcal{G}'$ are isomorphisms from the homology the graded module braid to graded module braid.
\end{definition}

Note that the internal $\I$, in the previous definition, may be $\PP$. Furthermore, we could have a collection of interval $\mathcal{I}\subseteq\I(<)$ such that for all $\I\in\mathcal{I}$ the previous diagram commutes, in this case, we say that $T$ \textit{weakly covers} $\theta$ \textit{on} $\mathcal{I}$. 

%\textcolor{blue}{It is worth to notice that when the boundary operator of both chain complex braids $\mathcal{C}$ and $\mathcal{C}'$ is trivial, it follows that there exists a transition matrix which cover an isomorphism $\theta$. Such result is used to guarantee the existence of classical topological transition matrix.}

\begin{theorem}\label{teo_weakly_cover}
Let $\Delta:\bigoplus_\PP C(p)\rightarrow \bigoplus_\PP C(p)$ and $\Delta':\bigoplus_\PP C'(p)\rightarrow \bigoplus_\PP C'(p)$ be connection matrices of $\mathcal{G}$ and $\mathcal{G}'$, respectively, and let $\theta:\mathcal{G}\rightarrow \mathcal{G'}$ be an isomorphism. If the order $<$ is stackable with a decomposition $\mathcal{I}=\{\I_1,...,\I_n\}$, then there is a transition matrix $T:\bigoplus_\PP C(p)\rightarrow \bigoplus_\PP C'(p)$ which weakly cover $\theta$ on $\mathcal{I}$.
\end{theorem}
\prooff It follows from the proof of theorem 3.5 in \cite{FM}, by checking that the induction process of applying Theorem 3.8 in \cite{FM} gives an collection of interval that equals to $\mathcal{I}$. \cqd

Note that trivial order and total order are stackable. Moreover, in \cite{FM} is given a definition of N-free order that helps to check if the order is stackable, in general, it is more likely a partial order be a stackable order than the opposite. However if $<$ is not stackable one still can analyze any interval $\I$ by applying Theorem \ref{teo_weakly_cover} in order to obtain a transition matrix $T(\I)$ that covers $\theta(\I)$, in other words, it is possible to get information between $p,q\in\PP$ by considering any interval which has $p$ and $q$.

Now we explore a special order $<_k$ that leads us to another existence result by using a similar idea coming from the Morse index of gradient-like flows. But before defining $<_k$, we present the following lemma. 

\begin{lemma}\label{Delta_block}
If there is at most one $k(p)$ such that $C_{k(p)} (p)\neq 0$, for all $p\in\PP$, then $\Delta(p,p')=0$ for all $p,p'\in\PP$ such that $|k(p)-k(p')|\neq 1$ or $p'<p$.
\end{lemma}
\prooff
Since $\Delta$ is a boundary map and $C(p)=C(p')=0$ except in dimension $k(p)$ and $k(p')$, we obtain that $\Delta(p,p')=0$ when $|k(p)-k(p')|\neq 1$. If $p'<p$, it follows by upper triangularity of $\Delta$ that $\Delta(p,p')=0$. \cqd

Under the hypothesis of Lemma \ref{Delta_block}, define $<_{k}$ a partial order in the set $\PP$ such that $$p<_k p'\Leftrightarrow k(p)<k(p').$$
The previous lemma seems to be quite obvious but plays the main role in the proof of the next theorem.

\begin{theorem}\label{almost_MS}
Let $\Delta:\bigoplus_\PP C(p)\rightarrow \bigoplus_\PP C(p)$ and $\Delta':\bigoplus_\PP C'(p)\rightarrow \bigoplus_\PP C'(p)$ be connection matrices of $\mathcal{G}(<_{k})$ and $\mathcal{G}'(<_{k})$, respectively, and let $\theta: \mathcal{G} \rightarrow \mathcal{G}'$ be an isomorphism between graded module braids. If there is at most one $k(p)$ such that $C_{k(p)} (p)\neq 0\neq C'_{k(p)}(p)$, for all $p\in\PP$, then there exists an unique transition matrix $T:\bigoplus_\PP C(p)\rightarrow \bigoplus_\PP C'(p)$ which covers $\theta$. 
\end{theorem}

\prooff
Lemma \ref{Delta_block} and the partial order $<_k$ guarantee that $T$ is a chain map. Moreover it allows us to use the proof of Theorem 6 in \cite{FdRV} to obtain the result.\cqd

Note that $T$ obtained in the previous theorem is a matrix in block form as in Theorem 6 in \cite{FdRV}. Furthermore, it is possible to consider a generic partial order $<$ instead of $<_k$ in Theorem \ref{almost_MS}, nevertheless the proof split in various cases where $<$ failed to be $<_k$. However, in case $<$, be aware that $T$ may not be unique neither a matrix in block form.

Observe that Theorem \ref{almost_MS} can be applied in a more general context than Theorem 6 in \cite{FdRV}. For instance, the flow on parameters $\lambda$ and $\mu$ do not need to be Morse-Smale without periodic orbits; we just need that the Morse decomposition has to meet the hypothesis that there is at most one $k(p)$ such that $CH_{k(p)} (M(p))\neq 0$, for all $p\in\PP$. Such a property is not satisfied for some Morse sets such as repeller periodic orbits and attractor periodic orbits, however saddle-saddle connections with same index satisfies this hypothesis.

\section{Algebraic Transition Matrix.}

The algebraic transition matrix theory is developed in \cite{FM}. There the authors define the algebraic transition matrix in the setting of a parameterized family of flows after having earlier developed an algebraic theory that they employ. Their algebraic transition matrix is an example of a similarity transformation defined between connection matrices. 

Let  $\Delta:\bigoplus_\PP C(p) \rightarrow \bigoplus C(p)$ be a connection matrix for a graded module braid $\mathcal{G}$ over a partial order $<$. If $T:\bigoplus C(p) \rightarrow \bigoplus C'(p)$ is $<$-upper triangular matrix such that $T(p)$ is an isomorphism for all $p\in\PP$, then it follows that $\Delta':= T\Delta T^{-1}$ is also a connection matrix for $\mathcal{G}$. In this way (as referred to above) $T$ can be thought of as a similarity transformation between connection matrices $\Delta$ and $\Delta'$.

\begin{proposition}\label{Tchain}\emph{[Franzosa]}
Let $C = \{C(p)\}_{p\in\PP}$ and $C' = \{C'(p)\}_{p\in\PP}$ be collections of
graded modules, and $\Delta:\bigoplus_\PP C(p)\rightarrow \bigoplus_\PP C(p)$ and $\Delta':\bigoplus_\PP C'(p)\rightarrow \bigoplus_\PP C'(p)$ be $<$-upper triangular boundary maps. If $T : \bigoplus_\PP C(p)\rightarrow \bigoplus_\PP C'(p)$ is $<$-upper triangular and such that $T\Delta=\Delta'T$, then $\{T(\I)\}_{\I\in \I(<)}$ is a chain map from $\mathcal{C}\Delta$ to $\mathcal{C}\Delta'$.
\end{proposition}

Here, for the purpose of unifying the various transition matrix theories--and being most general--we define
the above-mentioned similarity transformation between connection matrices as  follows:
\begin{definition}\label{ATM}
Let $\Delta:\bigoplus_\PP C(p)\rightarrow \bigoplus_\PP C(p)$ and $\Delta':\bigoplus_\PP C'(p)\rightarrow \bigoplus_\PP C'(p)$ be connection matrices.
If $T : \bigoplus_\PP C(p)\rightarrow \bigoplus_\PP C'(p)$ is $<$-upper triangular such that $T(p)$ is an isomorphism for all $p \in\PP$   and such that $\Delta' T=T\Delta$, then $T$ is called an \textbf{algebraic transition matrix} from $\Delta$ to $\Delta'$.

%Let $C = \{C(p)\}_{p\in(\PP,<)}$ and $C' = \{C'(p)\}_{p\in(\PP,<)}$ be collections of
%graded modules, if $T:\bigoplus_\PP C(p) \rightarrow \bigoplus_\PP C'(p)$ is upper triangular and such that $T(p)$ is an isomorphism for all $p\in\PP$, then $T$ is called \textbf{algebraic transition matrix}. 
\end{definition}
The algebraic transition matrix as defined in \cite{FM} {then} is an example of an algebraic transition matrix in Definition \ref{ATM} above. With the current definition of algebraic transition matrix, we then show that it is a transition matrix.

Note that, given $\Delta$, $\Delta'$ and $T$ as above, $T$ is a chain map between $C\Delta$ and $C\Delta'$ (see Proposition \ref{Tchain}). Furthermore $T$ induces homology isomorphisms $T_{\ast}(\I): H\Delta(\I) \rightarrow H\Delta'(\I)$ for all $\I \in \I(<)$. 
Now, if $\Delta$ is a connection matrix and $\Phi: H\Delta \rightarrow \mathcal{G}$ is a graded module braid isomorphism associated with $\Delta$, then if for each $\I \in \I(<)$ we define $\Phi'(\I)=\Phi(\I)\circ T_\ast^{-1}(\I)$, then the isomorphisms $\Phi'(\I)$ define a graded module braid isomorphism $\Phi'$ between $H\Delta'$ and $\mathcal{G}$. Also, clearly, for each $I \in \I(<)$, the following diagram commutes 
$$
\xymatrixcolsep{3pc} \xymatrix{
\mathcal{H}\Delta(\I) \ar[r]^{{T}_\ast(\I)} \ar[d]_{\Phi(\I)} & \mathcal{H}\Delta'(\I) \ar[d]^{\Phi'(\I)}\\
\mathcal{G}(\I) \ar[r]^{id(\I)}                   & \mathcal{G}(\I)
}
$$
where $id(\I)$ is the identity map on $\mathcal{G}(\I)$. Thus $T$ covers $id:\mathcal{G} \rightarrow \mathcal{G}$ (relative to $\Phi$ and $\Phi'$), and it follows that $T$ is a transition matrix.

It follows from above, that given a connection matrix $\Delta$ for a graded module braid $\mathcal{G}$, then other connection matrices for $\mathcal{G}$ can be obtained via transition matrices. The converse question is significant; that is, given connection matrices, $\Delta$ and $\Delta'$, is there a transition matrix for them? This question was addressed with some initial positive results in \cite{FM} and in the previous section \ref{Existence}. 

Define ATM$(<)$ to be the set of all algebraic transition matrices which are similarity transformations between connection matrices $\Delta:\bigoplus_{(\PP,<)} C(p) \rightarrow \bigoplus C(p)$ and $\Delta':\bigoplus_{(\PP,<)} C'(p) \rightarrow \bigoplus C'(p)$. From the existence theorem in section \ref{Existence}, one can prove the next theorem by applying the same idea from Theorem 4.3 in \cite{FM} and Lemma 3 in \cite{FdRV}.

\begin{theorem}\label{thrm_alg}
Let $M_\lambda=\{M_\lambda(\pi)\}_{\pi\in (\PP,<_\lambda)}$ and $M_\mu=\{M_\mu(\pi)\}_{\pi\in (\PP,<_\mu)}$ be Morse decompositions, $\Delta_\lambda$ and $\Delta_\mu$ their respective connection matrices. Moreover, assume that $M_\lambda$ and $M_\mu$ are related by continuation with an admissible ordering $<$.
If $<$ is stackable  or $<_k$ and $T_{\lambda,\mu}(\I_{j-1},\I_j)\neq0$ for all $T_{\lambda,\mu}\in$ ATM$(<)$ then there exists $s\in[0,1]$ such that $$C\left(M_{\omega(s)}(\I_{j}),M_{\omega(s)}(\I_{j-1})\right)\neq \emptyset,$$ where $\omega:[0,1]\rightarrow \Lambda$ is a path that continues $M_\lambda$ and $M_\mu$.\cqd
\end{theorem}

Notice in the previous theorem, when it is possible to set $\I_{j}$ and $\I_{j-1}$ to have only one element, Theorem \ref{thrm_alg} will appear more closely to the next theorems about bifurcation connections between Morse sets. Of course, the lack of covering the flow-defined continuation isomorphism $F$ imposes less features for algebraic transition matrices, since covering $F$ is intrinsically related how those connections occur along a path on $\Lambda$ as well as substantial existence results (see next sections).

\section{Singular Transition Matrix}

%This form of transition matrix is by far the most general. In fact, it is sufficient to assume the existence of an isolated invariant set that continues over ? and the existence of Morse decompositions of this invariant set at two different parameter values. 
In most cases singular transition matrices can only be computed via the dynamics of the slow system, since these transition matrices are essentially just submatrices of a connection matrix, see \cite{R}. Thus, in practice, one may find difficulty in obtaining those matrices if the objective is to understand the dynamics of the parametrized family. However, by showing that those matrices are transition matrices that cover an isomorphism, we actually can use the singular transition theory to assist in the development of other transition matrices, as we will show in the following sections.

Let 
\begin{equation}\label{eq1}
\dot x=f(x,\lambda),
\end{equation}
be a parametrized family of ordinary differential equations defined in $\mathbb{R}^n$, where the parameter space $\Lambda=\R$. Assume that the Morse decomposition $\mathcal{M}(S_\lambda)=\{M_\lambda(p)\in\PP\}$ continues over $\R$ and that connection matrices $\Delta_{-1}$ and $\Delta_1$ for the Morse decompositions $\mathcal{M}(S_{-1})$ and $\mathcal{M}(S_{1})$, respectively, are known. Moreover, let $N\subseteq \R^n$ be an isolating neighborhood for $S_\lambda$, $\lambda\in\R$.

One introduces slow dynamics in the parameter space of (\ref{eq1}) with the purpose of comprehending the bifurcations that occur for $-1<\lambda<1$. Hence (\ref{eq1}) can be written as:

\begin{equation}\label{eq_fss}
\begin{array}{c}
\dot x=f(x,\lambda),\\
\dot \lambda = \epsilon(\lambda^2-1)
\end{array}
\end{equation}

where $\epsilon>0$. Define
$$
M(p^+):=M_1(p),\ \ M(p^-):=M_{-1}(p)
$$
and $M(p^\pm):=M(p^+)\cup M(p^-)\subseteq \R^n\times\{\pm1\}$. For $\epsilon>0$ sufficiently small, $N\times [-2,2]$ is an isolating neighborhood for the flow $\phi_\epsilon$ generated by (\ref{eq_fss}). Let $K_\epsilon :=Inv(N\times[-2,2],\phi_\epsilon).$ Now observe that since $\dot\lambda<0$ if $\lambda\in(-1,1)$, for all $\epsilon>0$,
$$
\mathcal{M}(K_\epsilon)=\{M(p^\pm)\ |\ p\in\PP\}
$$
is a Morse decomposition, and there is an admissible ordering given by
$$
\begin{array}{c}
q^-<p^+,\\
q^-<p^-\ \Leftrightarrow \ q<_{-1}p,\\
q^+<p^+\ \Leftrightarrow \ q<_1 p,
\end{array}
$$
where $<_{-1}$ and $<_{1}$ are admissible orderings for $\mathcal{M}(S_{-1})$ and $\mathcal{M}(S_{1})$, respectively. Denoting a connection matrix for $\mathcal{M}(K_\epsilon)$ by $\Delta_\epsilon$, since the dynamics on the subspaces $\R^n\times \{\pm1\}$ are given exactly by the flows generated by $\dot x= f(x,\pm1)$, we have that
$$
\Delta_\epsilon: \bigoplus_{p\in\PP}H_\ast\left(M(p^-)\right) \bigoplus_{p\in\PP}H_\ast\left(M(p^+)\right) \rightarrow \bigoplus_{p\in\PP}H_{\ast-1}\left(M(p^-)\right) \bigoplus_{p\in\PP}H_{\ast-1}\left(M(p^+)\right)
$$
takes the form
\begin{equation}\label{delta_epsilon}
\Delta_\epsilon=
\left(
  \begin{array}{cc}
    \Delta_- & T_\epsilon \\
    0        & \Delta_+ \\
  \end{array}
\right)
\end{equation}
where $\Delta_-$ is the connection matrix for $\mathcal{M}(S_-)$ and $\Delta_+$ is the conjugation of degree 1 of the connection matrix of $\mathcal{M}(S_+)$ (this conjugation is necessary because on $\lambda=1$ we have an increase by one in the dimension of the unstable manifold, which introduces a suspension of the Conley index). %Observe that the space where $\Delta_\epsilon$ is defined is not a direct some of spaces where $\Delta_\pm$ are defined. 
A contribution of Reineck in \cite{R} was to formalize the expression (\ref{delta_epsilon}). When $\epsilon\rightarrow 0$ the limit $T_\epsilon$ is well defined and the resulting matrices are referred to as the \textbf{\emph{R-singular transition matrix}}. Ignoring the +1 conjugation on $\Delta_1$, it follows.

\begin{theorem}\emph{[Reineck]}
An R-singular transition matrix T from $\lambda=1$ to $\lambda=-1$ satisfies the following properties:
\begin{description}
  \item[(i)] $\Delta_{-1} T + T\Delta_1=0$;
  \item[(ii)] $T$ is an isormorphism;
  \item[(iii)] $T$ is an upper triangular matrix with respect to $<$;
  \item[(iv)] If $T(p,q):H(M(p))\rightarrow H(M(q))$ is nonzero, then there exists a finite sequence $1\geq\lambda_1\geq\lambda_2\geq\ldots \geq \lambda_k\geq0$ and corresponding $p_i\in \PP$ such that $p_i>_{\lambda_i} p_{i+1}$ where $>_{\lambda_i}$ is the flow defined order under $\phi_{\lambda_i}$. 
\end{description}
\end{theorem}

Now we define the singular transition matrix presented in \cite{MM2}. As one can note,  in \cite{R} and \cite{MM2}, different suspension isomorphisms to define singular transition matrix. This difference is important to note because as we point at in Remark \ref{remark} below, those matrices are transition matrices that cover different isomorphisms.

Following the same idea as in Reineck's development, Mischaikov and McCord in \cite{MM2} created a new parameter space that incorporates the drift flows and the one-parameter families in $\Lambda$. More specifically, let $\mathcal{D}^+=\mathcal{P}(\Lambda)  \times {\mathcal{G}}^+$ is a parameter space for flows on $X\times [-1,2]$, where $\mathcal{P}(\Lambda)=\{\alpha:[0,1]\rightarrow \mathbb{R}\}$ is the set of paths in $\Lambda$ and
$$
{\mathcal{G}}^{+}=\{g:[-1,2]\rightarrow \mathbb{R}\ |\ g\in C^0, g((-1,0)\cup (1,2))> 0> g((0,1)), g(-1)=g(0)=g(1)=g(2)=0, \text{ or } g\equiv 0\}.
$$

The choice of $[-1,2]$ instead of $[0,1]$ is convenient, since for the drift flow $\dot{s}=g(s)$ on $[0,1]$ has $0$ as a hyperbolic attractor and $1$ as a hyperbolic repeller, therefore agreeing with the drift flow in \cite{R}. Note that the domain of any path in $\mathcal{P}$ is $[0,1]$. We define 
$$
\tau(s)=\left\{\begin{array}{ccc}
-s & \text{for} & -1\leq s \leq 0,\\
s & \text{for} & 0\leq s\leq 1,\\
2-s & \text{for} & 1\leq s \leq 2.
\end{array}\right .
$$

For any $(\alpha,g)\in \mathcal{D}^+$ consider a flow on $X\times [-1,2]$ over $(\alpha,g)$ given by
$$
\begin{array}{ccl}
\dot{x}&=&f(x,\alpha \tau(s)),\\
\dot{s}&=&g(s).
\end{array}
$$
Observe that if $\alpha$ is a constant path $\lambda$, thus we obtain a product flow 
$$
\begin{array}{ccl}
\dot{x}&=&f(x,\lambda),\\
\dot{s}&=&g(s).
\end{array}
$$
When $||g||\rightarrow 0$, we have the original parameterized family of flows
$$
\dot{x}=f(x,\alpha\tau (s))
$$
restricted to the image of $\alpha$.

We are interested in studying the behavior on $\mathcal{D}^+([0,1])$, where the flow on $X\times [0,1]$ over $(\alpha,g)$ is the restriction of the flow on $X\times \Lambda$ to the one-parameter  family of flows picked out by $\alpha$. Let $C\subseteq [0,1]$ be a connected isolated set of zeros for $g\in \mathcal{G}^+$ restricted to $(-1,2)$, thus if $S$ is an isolated invariant set that continues over $\alpha[0,1]$ for some $\alpha$ then $S_\alpha=\{(x,s)\ |\ s\in C, x\in S_{\alpha(s)}\}$ is an isolated invariant set for the flow over $(\alpha,g)$. See \cite{MM2}.

Let $S_\alpha(C)$ denote the restriction of $S_\alpha$ to $C\subseteq[0,1]$, the following proposition gives us a way to compute Conley indices for isolated invariant sets over $\mathcal{D}^+$.

\begin{proposition}\emph{[Mischaikow-McCord]}
If $C$ is a connected, isolated set of zeros of $g$, then $h(C)$, the homotopy Conley index of $C$ in $[-1,2]$, is either $\Sigma^1, \Sigma^0$ or $\bar{0}$. Let $\mathcal{C}_g$ denote the component of $\mathcal{G}^+$ that contains $g$. If $S$ is an isolated invariant set that continues over $\Lambda$ then, over $\mathcal{P}\times \mathcal{C}_g$, the Conley index of $S_\alpha(C)$ in $X\times[-1,2]$ is $h(S)\wedge h(C)$.
\end{proposition}

By the previous proposition, we have that the homology Conley index of $S_\alpha(C)$ is the tensor product
$$
CH(S_\alpha(C))=CH(S)\otimes CH(C).
$$
In this setting, we present the \textit{index suspension isomorphism} $\Sigma(S)$ defined in \cite{MM2} by the following composition
$$
CH_k(S_{\alpha(c)})\xrightarrow{\otimes \sigma_n} CH_k(S_{\alpha(c)})\otimes CH_n(C)\xrightarrow{\times} CH_{n+k}(S_{\alpha(c)}\times C)=CH_{n+k}(S_{\bar{\alpha}(c)}(C))\xrightarrow{F_{\alpha,\bar{\alpha}(c)}}CH_{k+n}(S_\alpha(C)),
$$
where $c\in C$, $\alpha$ is a path, $\bar{\alpha}(c)$ is the constant path $\alpha(c), \sigma_n$ is the generator of $CH_n(C)$ and $F_{\alpha,\bar{\alpha}(c)}$ is the continuation isomorphism along a path in $\mathcal{D}^+$ from $(\bar{\alpha}(c),g)$ to $(\alpha,g)$. More concisely, we have
$$
\Sigma(S):CH_k(S_{\alpha(c)})\xrightarrow{F_{\alpha,\bar{\alpha}(c)}\circ \times\circ\otimes\sigma_n} CH_{k+n}(S_\alpha(C)).
$$
Given that $M(\textbf{I})$ is an isolated invariant set which also continues, we can define 
$$\Sigma(\I)=\Sigma(M(\textbf{I})):CH_k(M(\I)_{\alpha(c)})\xrightarrow{F_{\alpha,\bar{\alpha}(c)}\circ \times\circ\otimes\sigma_n} CH_{k+n}(M(\I)_\alpha(C)).$$

From the dynamics in $[0,1]$, note that, $S_{(\alpha,g)}$ has $S_{\alpha(0)}$ as an attractor and $S_{\alpha(1)}$ as a repeller. These characteristics will be the same for all $g\in\mathcal{G}^+$, even though the structure of of the connecting orbit set may vary with $g$.  Computing the connection matrix for the flow in $X\times\alpha[0,1]$, we have
$$
\Delta_g: \bigoplus_{p\in\PP}H_\ast\left(M(p)_{\alpha(0)}\right) \bigoplus_{p\in\PP}H_\ast\left(M(p)_{\alpha(1)}\right) \rightarrow \bigoplus_{p\in\PP}H_\ast\left(M(p)_{\alpha(0)}\right) \bigoplus_{p\in\PP}H_\ast\left(M(p)_{\alpha(1)}\right)
$$
takes the form
\begin{equation}\label{delta_g}
\Delta_g=
\left(
  \begin{array}{cc}
    \Delta_{\alpha(0)} & T_g \\
    0        & \Delta^\Sigma_{\alpha(1)} \\
  \end{array}
\right)
\end{equation}

where $\Delta_{\alpha(0)}$ is a connection matrix for $\mathcal{M}(S)_{\alpha(0)}$ and $\Delta^\Sigma_{\alpha(1)}$ is the conjugation by $\Sigma$ of a connection matrix of $\mathcal{M}(S)_{\alpha(1)}$, {see Lemma \ref{last_dia}}. 

When    $||g||\rightarrow 0$ the limit $T_g$ gives us a {\bf{MM-singular transition matrix}} $T_s$. Observe that $\Delta_g$ may differ from $\Delta_\epsilon$ since we used different suspension isomorphisms to define these singular transition matrices.

%Letting $||g||\rightarrow 0$, we can do the same process, needed to define singular transition matrix, to obtain another singular transition matrix 
%$\Delta_g$ has a different structure which $\Delta_\epsilon$ may not have, it is the following property.
The following theorem points at an important difference between the structure of $\Delta_g$ and the structure of $\Delta_\epsilon$.

\begin{theorem}\emph{[McCord-Mischaikow]}\label{comm_F}
The connection homomorphism for the attractor-repeller decomposition $\left(\mathcal{M}_{\alpha(0)}(\I),\mathcal{M}_{\alpha(1)}(\I)\right)$ of $\mathcal{M}_{(\alpha,g)}(\I)$ is an isomorphism, that is computed by continuation of $\mathcal{M}_{(\alpha,g)}(\I)$ across $\Lambda$. That is, there is a commutative diagram
$$
\xymatrixcolsep{3pc}\xymatrix{
H_{k,\alpha(1)} (\I)    \ar[d]^{id}       \ar[r]^{\Sigma(\I)}                &  H_{k+1,\alpha(1)} (\I)       \ar[r]^{\delta}                 & H_{k,\alpha(0)} (\I) \ar[d]^{id} \\
H_{k,\alpha(1)} (\I)                                    \ar[rr]^{F_{\alpha(1),\alpha(0)}(\I)}             &            &               H_{k,\alpha(0)} (\I)
}
$$
\end{theorem}

As a consequence we have that the flow defined map of attractor-repeller map $\delta$ is an isomorphism. The following theorem addresses both the R-singular and the MM-singular transition matrix.

\begin{theorem}
A singular transition matrix is a 1-degree transition matrix. % which covers the flow-defined continuation isomorphism $F$.
\end{theorem}

\begin{remark}\label{remark}
As we can see in the following proof, Theorem \ref{STM=TM} works for both R-singular and the MM-singular transition matrix.  More specifically, they are 1-degree transition matrices that cover the 1-degree isomorphisms $\delta_\Psi$ and $F\circ\Sigma^{-1}$, respectively.
%, where $\delta$ is the flow defined map of attractor-repeller pair after the repeller Conley index at parameter $\lambda$ being suspended by the suspension isomorphism $\Psi$ coming from Mayer-Vietoris sequence and $\Sigma$ is the isomorphism coming from the composition of $F$ and the suspension isomorphism for a constant path in $\Lambda$.
\end{remark}

\prooff
In order to simplify notation, we define $\lambda=\alpha(0)$ and $\mu=\alpha(1)$. Let $T_s^1$ and $T_s$ be the R-singular and MM-singular transition matrices {(defined by $\Psi$ and $\Sigma$)}, respectively. Choose an interval $\I\in\mathcal{I}(<)$, where $<$ is the flow defined order for the product flow. Hence ($\I_\mu, \I_\lambda$) is an attractor-repeller pair and, by $\Delta_g$ being a connection matrix, we have that the following diagram commutes

$$
\xymatrixcolsep{5pc}\xymatrix{
 H_{k+1}\Delta_\lambda^\Sigma (\I)  \ar[d]^{\Phi_\lambda^\Sigma(\I)}
 \ar[r]^{[T_{s}]} & H_k\Delta_\mu (\I) \ar[d]^{ \Phi_\mu (\I)} \\
  H_{k+1,\lambda} (\I)       \ar[r]^{\delta}                 & H_{k,\mu} (\I).
}
$$

Therefore $T_s$ is a 1-degree transition matrix that covers $\delta=F\circ\Sigma^{-1}$ since $T_s$ is a 1-degree chain map and $\delta$ is defined by a 1-degree graded module braid isomorphism. If one changes $\Sigma$ for $\Psi$ we obtain $T^1_s$ is a 1-degree transition matrix covering $\delta_\Psi$. Note that in this case we may not have $\delta_\Psi=F\circ\Psi^{-1}$ since $\delta_\Psi$ is a flow defined map of the attractor-repeller pair for which the Conley index of the repeller was suspended by $\Psi$.\cqd

%Theorem \ref{STM=TM} works for both definition of singular transition matrix found in \cite{R1} and \cite{MM2}, as one can see in the its proof.  More especifically, STM from \cite{R1} and \cite{MM2} are transition matrices which covers $\delta_\Phi$ and $F\circ\Sigma^{-1}$ respectively.

\begin{lemma}\label{last_dia}
$\Sigma:\mathcal{H}_{\lambda,\ast}\rightarrow \mathcal{H}_{\lambda,\ast+1}$ is a braid isomorphism between Conley index braids and $\Delta^\Sigma_\lambda$ is $\Sigma$ conjugated to a connection matrix $\Delta_\lambda$ at parameter $\lambda$, in other words the following diagram commutes
$$
\xymatrixcolsep{5pc}\xymatrix{
H_k\Delta_\lambda (\I)  \ar[d]^{\Phi_\lambda (\I)}    \ar[r]^{[\oplus_{\pi\in\I}\Sigma(\pi)]}   & H_{k+1}\Delta_\lambda^\Sigma (\I)  \ar[d]^{\Phi_\lambda^\Sigma(\I)}\\
H_{k,\lambda} (\I)        \ar[r]^{\Sigma(\I)}                &  H_{k+1,\lambda} (\I)         
}
$$
\end{lemma}
\prooff

$\Sigma$ is an isomorphism by definition. Now it reminds to prove that for all adjacent pair $(\I,\J)$ the following diagram commutes
$$
\xymatrixcolsep{2pc}\xymatrix{
\cdots \ar[r] & H_{\lambda,\ast}(\I) \ar[r] \ar[d]^{\Sigma(\I)} & H_{\lambda,\ast}(\I\J) \ar[r] \ar[d]^{\Sigma(\I\J)} & H_{\lambda,\ast}(\J) \ar[r] \ar[d]^{\Sigma(\J)} & \cdots  \\ 
\cdots \ar[r] & H_{\lambda,\ast+1}(\I) \ar[r] & H_{\lambda,\ast+1}(\I\J) \ar[r] & H_{\lambda,\ast+1}(\J) \ar[r] & \cdots
}
$$
Indeed, the commutativity follows from the same property that the flow-defined Conley index isomorphism $F$ has. In other words, $\Sigma$ is defined via $F$, thus one obtains following diagram
$$
\xymatrixcolsep{2pc}\xymatrix{
\cdots \ar[r] & H_{\lambda,\ast}(\I) \ar[r] \ar[d]^{\times\circ\otimes\sigma(\I)} & H_{\lambda,\ast}(\I\J) \ar[r] \ar[d]^{\times\circ\otimes\sigma(\I\J)} & H_{\lambda,\ast}(\J) \ar[r] \ar[d]^{\times\circ\otimes\sigma(\J)} & \cdots  \\ 
\cdots \ar[r] & H_{\lambda,\ast+1}(\I) \ar[r] \ar[d]^{F(\I)} & H_{\lambda,\ast+1}(\I\J) \ar[r] \ar[d]^{F(\I\J)} & H_{\lambda,\ast+1}(\J) \ar[r] \ar[d]^{F(\J)} & \cdots  \\ 
\cdots \ar[r] & H_{\lambda,\ast+1}(\I) \ar[r] & H_{\lambda,\ast+1}(\I\J) \ar[r] & H_{\lambda,\ast+1}(\J) \ar[r] & \cdots
}
$$
which commutes since $\times\circ\otimes\sigma$ and $F$ are isomorphism between braids.  Therefore $\Sigma$ is a graded module braid isomorphism.

Note that, in order to make a cleaned proof, we omitted some subscripts which come from the definition of $\Sigma$. But be aware that for a small enough neighborhood the Conley index $H_{\lambda,\ast}(\I)$ was suspended by $\times\circ\otimes\sigma$ for a constant path and soon after the suspended Conley index $H_{\lambda,\ast+1}(\I)$ was homotoped via $F$ to a path $\alpha$ inside a small enough neighborhood.

Defining $$\Delta_{\lambda}(\pi',\pi):=\Sigma^{-1}(\pi)\circ\Delta^\Sigma(\pi',\pi)\circ\Sigma(\pi')
$$
we obtain that $\Delta_\lambda$ is a connection matrix at parameter $\lambda$ since $\Sigma:\mathcal{H}_{\lambda,\ast}\rightarrow \mathcal{H}_{\lambda,\ast+1}$ is an isomorphism between Conley index braids, therefore the diagram in Lemma \ref{last_dia} commutes. \cqd

MM-Singular transition matrix $T_s$ defined via $\Sigma$ has another feature as follows by the next theorem.

\begin{theorem}\label{STM=TM}
The MM-singular transition matrix $T_s$ composed with the induced map of $\bigoplus_{\pi\in\PP} \Sigma(\pi)$ is a transition matrix
$$T=T_s\circ \left(\bigoplus_{\pi\in\PP} \Sigma(\pi)\right)$$ that covers Conley flow-defined isomorphim $F$.
\end{theorem}
\prooff

Suppose that the following diagram commutes for every interval $\I\in\mathcal{I}(<)$

$$
\xymatrixcolsep{5pc}\xymatrix{
H_k\Delta_\lambda (\I)  \ar[d]^{\Phi_\lambda (\I)}    \ar[r]^{[\oplus_{p\in\I}\Sigma_p]}   & H_{k+1}\Delta_\lambda^\Sigma (\I)  \ar[d]^{\Phi_\lambda^\Sigma(\I)}
 \ar[r]^{[T_{s}]} & H_k\Delta_\mu (\I) \ar[d]^{ \Phi_\mu (\I)} \\
H_{k,\lambda} (\I)    \ar[d]^{id}       \ar[r]^{\Sigma(\I)}                &  H_{k+1,\lambda} (\I)       \ar[r]^{\delta}                 & H_{k,\mu} (\I) \ar[d]^{id} \\
H_{k,\lambda} (\I)                                    \ar[rr]^{F_{\lambda,\mu}(\I)}             &            &               H_{k,\mu} (\I)
}
$$
\begin{center}
Diagram 1
\end{center}
%\begin{table}
%\caption{}\label{dia_big}
%\end{table}
where $[\oplus_{\pi\in\I}\Sigma(\pi)]$ and $[T_s]$ are induced isomorphisms of the chain maps $\oplus_{\pi\in\I}\Sigma(\pi)$ and $T_s$ respectively, and $\Phi$'s are isomorphisms which came from the definition of connection matrix. See Lemma \ref{last_dia} for a precise definition of $\Delta_\lambda$ and $\Phi_\lambda$.

Thus we have that the chain map $T=T_s\circ  \oplus_{\pi\in\I}\Sigma(\pi)$ covers the flow-defined Conley-index isomorphism $F_{\lambda,\mu}$, therefore $T$ is an GTTM and $T_s$ is a GTM which covers $F_{\lambda,\mu}\circ\Sigma^{-1}=\delta$.

Now we just need to prove that the latter diagram commutes. Indeed, 
$$
\xymatrixcolsep{5pc}\xymatrix{
H_{k,\lambda} (\I)    \ar[d]^{id}       \ar[r]^{\Sigma(\I)}                &  H_{k+1,\lambda} (\I)       \ar[r]^{\delta}                 & H_{k,\mu} (\I) \ar[d]^{id} \\
H_{k,\lambda} (\I)                                    \ar[rr]^{F_{\lambda,\mu}(\I)}             &            &               H_{k,\mu} (\I)
}
$$
commutes by Theorem \ref {comm_F} and
$$
\xymatrixcolsep{5pc}\xymatrix{
 H_{k+1}\Delta_\lambda^\Sigma (\I)  \ar[d]^{\Phi_\lambda^\Sigma(\I)}
 \ar[r]^{[T_{s}]} & H_k\Delta_\mu (\I) \ar[d]^{ \Phi_\mu (\I)} \\
  H_{k+1,\lambda} (\I)       \ar[r]^{\delta}                 & H_{k,\mu} (\I)
}
$$
commutes by definition of $\Delta_g$. The last diagram commutes by Lemma \ref{last_dia}.

Now one needs to check that the induced map $[T]=\left[T_s\circ \bigoplus_{\pi\in\PP} \Sigma(\pi)\right]$ is a well defined map between the braids $\mathcal{H}\Delta_\lambda$ and $\mathcal{H}\Delta_\mu$, indeed this follows since all maps from Diagram 1 are graded module braid maps. \cqd

Note that $\Sigma$ depends on the flow defined isomorphism $F$, therefore $\Delta^\Sigma_\lambda$ and $T_s$ depends on $F$ too, if one is willing to remove that dependence (in order to obtain those maps in an easy way) the next corollary solve this problem.

\begin{corollary}\label{indep}
The suspension isomorphism $\Sigma$ does not depend on $F$ when the flow at parameter $\lambda$ is structural stable. Moreover,
$$T=T_s\circ \otimes \sigma,
$$
where $T$ is a transition matrix that covers the Conley flow-defined isomorphism $F$ and $T_s$ is an MM-singular transition matrix.
\end{corollary}

\prooff
Let $\alpha$ be a path between $\lambda$ to $\mu$. If $\lambda=\mu$, let $\alpha$ be a constant path, hence we have that $\Sigma=\times \circ \otimes \sigma$, therefore $\Sigma$ does not depend on $F$ and by Theorem \ref{STM=TM} the result follows. 
Suppose that $\alpha$ is a non-constant path, since the flow in parameter $\lambda$ is structural stable there is a $\epsilon>0$ sufficient small such that $\alpha[0,\epsilon]$ does not have bifurcations, thus $F_{\overline{\alpha}(0)[0,\epsilon]}=id$, where $\overline{\alpha}(0)$ is a constant path evaluated in $\alpha(0)$.
 
By $\Lambda$ being simply connected, we have that there exists a homotopy $A:\alpha \simeq \xi$, where $\xi=\overline{\alpha}(0)[0,\epsilon]\ast\gamma\ast\alpha[\epsilon,1]$ and $\gamma$ is a path between $\alpha(0)$ and $\alpha(\epsilon)$. Hence $F_{\alpha[0,\epsilon]}=id$, thus $\Sigma$ does not depend on $F$, since we can suspended Conley index by $\Sigma$ locally. Therefore, by Theorem \ref{STM=TM}, the result follows. \cqd

Theorem \ref{strucStable} in the next section shows us the importance when $\Sigma$ does not depend on $F$, i. e., when one can apply Corollary \ref{indep}.

%Let $\Xi$ be the isomorphism, which comes from the definition of $\Delta_g$. For an interval $\I_\mu$ we have that $\Xi(\I_\mu)=\Phi_\mu(\I_\mu)$, since there was no suspension of Conley index at parameter $\mu$ (the equality is an abuse of notation because the $\Xi(\I_\mu)$ is actually $\Phi_\mu(\I_\mu)$ conjugated by induced inclusion map).

\section{Topological Transition Matrix.}

In this section, we develop more proprieties of generalized topological transition matrices from \cite{FdRV} given that such matrices are transition matrix which covers flow-defined Conley-index isomorphisms. Therefore we can use the interrelationship between singular and topological transition matrix theory in order to obtain refined properties, for instance, Theorem \ref{teo_prop1} is an improvement of Theorem 5 in \cite{FdRV} by no requiring existence assumption. Moreover, Theorem \ref{strucStable} describe appropriately the importance of the unify theory of transition matrix.

Let $M_{\lambda}=\{M_{\lambda}(\pi)\}_{\pi\in\mathbf{P}}$ and $M_{\mu}=\{M_{\mu}(\pi)\}_{\pi\in\mathbf{P}}$ be Morse decompositions, related by continuation, for the isolated invariant sets $S_\lambda\subseteq X_\lambda$ and $S_\mu\subseteq X_\mu$, respectively. 

\begin{definition}\label{GTTM}
If $T$ is a transition matrix that covers the Conley flow-defined isomorphism $F$, then we refer to $T$ as a \textbf{(generalized) topological transition matrix}.
\end{definition}

In order to simplify one may omit the term in parenthesis ``generalized" from Definition \ref{GTTM}, nevertheless be aware that this definition is an extension of topological transition matrix defined in \cite{MM1}, where they assume that there are no connections at the initial and final parameters of a continuation.

By the definition of transition matrix we obtain the following statement.

\begin{theorem}
Generalized topological transition matrix is a transition matrix which covers the flow-defined isomorphism $F$. \hfill $\square$
\end{theorem}

Here is another way to characterize generalized topological transition matrices.

\begin{proposition}\label{def_new}
$T$ is a \textbf{generalized topological transition matrix} related to the connection matrices $(\Delta_\lambda,\ \Phi_\lambda)$ and $(\Delta_\mu,\ \Phi_\mu)$, if and only if
$$T:\displaystyle \bigoplus_{p\in \PP}CH_\ast (M_\lambda(p)) \rightarrow \displaystyle \bigoplus_{ p\in  \PP}CH_\ast (M_\mu( p))$$
is a zero degree map such that

\begin{itemize}
  \item $\{T(\I)\}_{\I\in\I(<)}$ is a chain map from $\mathcal{C}\Delta_\lambda$ to $\mathcal{C}\Delta_\mu$;
  \item the following diagram commutes
\end{itemize}
\begin{table}[h!]\centering
\begin{tikzpicture}
  \matrix (m) [matrix of math nodes, row sep=3em,
    column sep=0.5em]{
	& H\Delta_\lambda(\emph{\I}) & & H\Delta_\lambda(\emph{\I\J}) & & H\Delta_\lambda(\emph{\J}) & & H\Delta_\lambda(\emph{\I})\\
	 H_\lambda(\emph{\I})       &  & H_\lambda(\emph{\I\J})       & & H_\lambda(\emph{\J})       & & H_\lambda(\emph{\I})     & \\
	& H\Delta_\mu(\emph{\I}) & & H\Delta_\mu(\emph{\I\J}) & & H\Delta_\mu(\emph{\J}) & & H\Delta_\mu(\emph{\I}) \\
	 H_\mu(\emph{\I})       & & H_\mu(\emph{\I\J})       & & H_\mu(\emph{\J})       & & H_\mu(\emph{\I})      &\\};	
  \path[-stealth]
    (m-1-2) edge (m-1-4) edge  (m-2-1)
            edge [densely dotted] (m-3-2)
    (m-2-1) edge [-,line width=6pt,draw=white] (m-2-3)
            edge (m-2-3) edge node [left] {{\small $F_{\lambda\mu}(\emph{\I})$}} (m-4-1)
    (m-3-2) edge [densely dotted] (m-3-4)
            edge [densely dotted] (m-4-1)
    (m-4-1) edge (m-4-3)
    (m-1-4) edge [densely dotted] (m-3-4) edge (m-2-3) edge (m-1-6)
    (m-3-4) edge [densely dotted] (m-4-3) edge [densely dotted] (m-3-6)
    (m-2-3) edge [-,line width=6pt,draw=white] (m-2-5) edge [-,line width=3pt,draw=white] (m-4-3)
            edge (m-4-3) edge (m-2-5)
	(m-4-3) edge (m-4-5)
	(m-1-6) edge (m-2-5) edge [densely dotted] (m-3-6) edge node [above] {{\small $\Delta_\lambda (\emph{\J},\emph{\I})$}} (m-1-8)
	(m-2-5) edge [-,line width=30pt,draw=white] (m-2-7) edge [-,line width=3pt,draw=white] (m-4-5) edge (m-4-5) edge node [above] {{\small $\delta_\lambda (\emph{\J},\emph{\I})$}} (m-2-7)
	(m-3-6) edge [densely dotted] (m-4-5) edge [densely dotted] (m-3-8)
	(m-4-5) edge (m-4-7)
	(m-1-8) edge (m-2-7) edge node [right] {{\small $\hat T(\emph{\I})$}} (m-3-8)
	(m-2-7) edge [-,line width=3pt,draw=white] (m-4-7) edge (m-4-7)
	(m-3-8) edge (m-4-7);
\end{tikzpicture}
\caption{}\label{dia_prin}
\end{table}
for all adjacent pairs \emph{$(\I,\J)$}, where $\hat T(\cdot)$ is the induced homology map of $T(\cdot)$.
\end{proposition}

Denote GTTM$(<)$ as the set of all generalized topological transition matrices with the partial order $<$.

When there are no connections in the $\lambda$ and $\mu$ parameters, then $\Delta_\lambda =0 =\Delta_\mu$. By Conley's theory we have that there is an isomorphism $\Phi_\lambda:C_\ast\Delta_\lambda(\PP) \rightarrow H_{\ast,\lambda}(\PP)$ for $\lambda\in\Lambda'$, where $C_\ast\Delta_\lambda(\PP)=\bigoplus_{\pi\in\PP}CH(M_\lambda(\pi))$ is the chain complex with connection matrix $\Delta_\lambda$.

Therefore, we can carry out the continuation along the path $\omega$ in two ways: first by continuing $S_\lambda$ along the path $\omega$ using the isomorphism $F_{\lambda,\mu}$; secondly continuing $\bigcup_{p\in \PP} M_\lambda(p)$ along the path $\omega$ by using isomorphism $E_{\lambda,\mu}=\bigoplus_{p\in \PP} F_{\lambda,\mu}(M(p))$. More precisely, we have the following diagram
$$
\xymatrixcolsep{3pc}\xymatrix{
C\Delta_\lambda(\PP) \ar[r]^{E_{\lambda,\mu}} \ar[d]^{\Phi_\lambda} &  C\Delta_\mu(\PP) \ar[d]^{\Phi_\mu}\\
H_{\lambda}(\PP) \ar[r]^{F_{\lambda,\mu}}  &  H_{\mu}(\PP)
}
$$

In general the diagram above is not commutative, in order to commute just define the following generalized topological transition matrix $T_{\lambda,\mu}=\Phi_\mu^{-1}\circ F_{\lambda,\mu}\circ \Phi_\lambda$, then we obtain from Proposition \ref{def_new} the following diagram commutes
$$
\xymatrixcolsep{3pc}\xymatrix{
C\Delta_\lambda(\PP) \ar[r]^{T_{\lambda,\mu}(\PP)} \ar[d]^{\Phi_\lambda} &  C\Delta_\mu(\PP) \ar[d]^{\Phi_\mu}\\
H_{\lambda}(\PP) \ar[r]^{F_{\lambda,\mu}(\PP)}  &  H_{\mu}(\PP)
}
$$
In this particular case, such matrix $T_{\lambda,\mu}$ is called a \textbf{(\textit{classical}) topological transition matrix}.

%The following theorem in \cite{MM1} summarizes some important properties for this matrix, which we refer to, from now on, as the classical topological transition matrix.
%
%\begin{theorem}\label{teo_classical_top}\emph{[McCord-Mischaikow]} Let $\Lambda'\subseteq \Lambda$ be such that for all $\lambda$ and $\mu\in\Lambda'$ there are no connection orbits in $M_\lambda$ and $M_\mu$, and $M_\lambda$ and $M_\mu$ are related by continuation. Then
%\begin{description}
%  \item[(i)] $\Delta_\mu T_{\lambda,\mu} + T_{\lambda,\mu}\Delta_\lambda = 0$;
%  \item[(ii)] $T_{\lambda,\mu}$ is an isomorphism;
%  \item[(iii)] $T_{\lambda,\mu}$ is upper triangular matrix with respect to order $<$ ;
%  \item[(iv)] If $\nu\in\Lambda'$ then $T_{\lambda,\lambda}=id$, $T_{\lambda,\nu}=T_{\mu,\nu}\circ T_{\lambda,\mu}$ and $T_{\mu,\lambda}=T_{\lambda,\mu}^{-1}$;
%  \item [(v)] If $T_{\lambda,\mu}(p,q)\neq0$ and $\omega$ is a path between $\lambda$ and $\mu$, then there exists a finite sequence $0<s_1\leq s_2\leq \ldots \leq s_{n}<1$ and a sequence $(p_i)\subseteq\PP$ such that $p_0=q,\ p_{n}=p$ and the connecting orbit set $C\left(M_{\omega(s_i)}(p_{i-1}),M_{\omega(s_i)}(p_i)\right)$ is nonempty.
%\end{description}
%\end{theorem}

Applying Theorem \ref{STM=TM} we have the following existence result.

\begin{theorem}\label{existence}
Let $M_\lambda=\{M_\lambda(\pi)\}_{\pi\in (\PP,<_\lambda)}$ and $M_\mu=\{M_\mu(\pi)\}_{\pi\in (\PP,<_\mu)}$ be Morse decompositions, assume that $M_\lambda$ and $M_\mu$ are related by continuation with an admissible ordering $<$. Then there exists a generalized topological transition matrix $T$ related to a $\Delta_\lambda\in\mathcal{CM}(<_\lambda)$ and a $\Delta_\mu\in\mathcal{CM}(<_\mu)$. \hfill $\square$
\end{theorem}
Now we can remove the existence hypotheses from item (v) of Theorem 5 in \cite{FdRV} using Theorem \ref{existence}.

\begin{theorem}\label{teo_prop1}
Let $M_\lambda=\{M_\lambda(\pi)\}_{\pi\in (\PP,<_\lambda)}$ and $M_\mu=\{M_\mu(\pi)\}_{\pi\in (\PP,<_\mu)}$ be Morse decompositions, $\Delta_\lambda$ and $\Delta_\mu$ their respective connection matrices. Moreover, assume that $M_\lambda$ and $M_\mu$ are related by continuation with an admissible ordering $<$. Then the generalized topological transition matrix $T$ satisfies the following properties:
\begin{description}
  \item[(i)] $T \circ \Delta_\lambda=\Delta_\mu\circ  T$;
  \item[(ii)] $T_{\lambda,\mu}(\{p\})=id$ and $T$ is upper triangular with respect to $<$;
  \item[(iii)] $T$ is an isomorphism;
  \item[(iv)] $T_{\lambda,\lambda}=id$, $T_{\lambda,\nu}(\I)=T_{\mu,\nu}\circ T_{\lambda,\mu}(\I)$ and $T_{\mu,\lambda}(\I)=T_{\lambda,\mu}^{-1}(\I)$ are generalized topological transition matrices, for all intervals $\I\in \cal{I}$ and $p\in\PP$, in particular $T=T(\PP)$.
  \item [(v)] Let $\omega:[0,1]\rightarrow \Lambda$ be a path that continues $M_\lambda$ to $M_\mu$. Assume that $T_{\lambda,\mu}(p,q)\neq0$ for all generalized topological transition matrices. Then there exists a finite sequence $0\leq s_1\leq s_2\leq \ldots \leq s_{n}\leq 1$ and a sequence $(p_i)\subseteq\PP$ such that $p_0=q,\ p_{n}=p$ and the set of connecting orbits $C\left(M_{\omega(s_i)}(p_{i-1}),M_{\omega(s_i)}(p_i)\right)$ is non-empty.
\end{description}
\end{theorem}
\prooff
From Theorem \ref{existence} we have the existence of $T_{\lambda,\mu}$ related to $\bar{\Delta}_\lambda\in\mathcal{CM}(<_\lambda)$ and $\bar{\Delta}_\mu\in\mathcal{CM}(<_\mu)$. However, $\Delta_\alpha$ may be different from $\bar{\Delta}_\alpha$, with $\alpha=\{\lambda, \mu\}$, since Theorem \ref{existence} does not specify which connection matrix $T$ is related. Therefore, one needs to check in the proof of Theorem 5 in \cite{FdRV} if all connection matrices at parameters $\omega(0),\ \omega(s_{\xi_1}),\ \ldots,\ \omega(1)$ agree.

Without loss of generality, assume that we are in case where the connection matrices at parameters $\omega(0)=\lambda,\ \omega(s_{\xi_1})=\nu,\ \omega(1)=\mu$ must agree. Indeed,  by Theorem \ref{existence} one obtains $T_{\lambda,\nu}$ related to $\Delta'_\lambda$ and $\Delta'_\nu$ and $T_{\nu,\mu}$ related to $\Delta''_\nu$ and $\Delta'_\mu$. Although, $\Delta'_\nu$ may be different to $\Delta''_\nu$, one can obtain the same matrix by using Theorem 4.8 in \cite{Fr2}. In other words, such theorem stats the existence of connection matrix by doing an induction process, since Theorem \ref{existence} use this result to obtain $T$ then one can do the same induction process starting from $\Delta'_\nu$ to obtain 
\begin{equation}
\Delta=
\left(
  \begin{array}{cc}
    \Delta'_\nu & T'_{\nu,\mu} \\
    0        & \Delta''_\mu\\
  \end{array}
\right),
\end{equation} where $\Delta''_\mu$ is a connection matrix at parameter $\mu$ and $T'_{\nu,\mu}$ is a generalized topological transition matrix related to $\Delta'_\nu$ and $\Delta''_\mu$. 

Note that, $T_{\nu,\mu}$ and $\Delta'_\mu$ may be different from $T'_{\nu,\mu}$ and $\Delta''_\mu$, respectively. But it is not a problem because one can obtain a $T'_{\lambda,\mu}=T_{\lambda,\nu}\circ T'_{\nu,\mu}$ and by hypothesis $T'_{\lambda,\mu}(p,q)\neq0$.\cqd

\begin{remark}
Applying the same idea, in the previous proof, on Theorem \ref{existence} and on equation \ref{delta_epsilon} and \ref{delta_g}, one can freely choose a connection matrix (in their respective set of connection matrices) to be $\Delta_\lambda$, $\Delta_-$ and $\Delta_{\alpha(0)}$ respectively.  However, $\Delta_\mu$, $\Delta_+$ and $\Delta^\Sigma_{\alpha(1)}$ are not free to choice, in fact they come from the inductive process of Theorem 4.8 in \cite{Fr2}. In other words, existence results coming from the singular transition matrix theory give us liberty to choose a connection matrix at the parameter that was not suspended, yet we can not choose freely a connection matrix at the parameter which was suspended.
\end{remark}

Using both topological and singular transition matrix, one is able to obtain a richer dynamical information in Theorem \ref{teo_prop1} without assuming that  $T_{\lambda,\mu}(p,q)\neq0$ for all generalized topological transition matrices. In other words, it is enough to assume $T_{\lambda,\mu}(p,q)\neq0$ for at least one matrix in GTTM.

\begin{theorem}\label{strucStable}
Let $M_\lambda=\{M_\lambda(\pi)\}_{\pi\in (\PP,<_\lambda)}$ and $M_\mu=\{M_\mu(\pi)\}_{\pi\in (\PP,<_\mu)}$ be Morse decompositions related by continuation with an admissible ordering $<$. Moreover, assume that the flow at parameter $\lambda$ is structural stable and a generalized topological transition matrix $T$ has a nonzero entrance $T_{\lambda,\mu}(p,q)\neq0$. Then there exists a finite sequence $0< s_1\leq s_2\leq \ldots \leq s_{n}\leq 1$ and a sequence $(p_i)\subseteq\PP$ such that $p_0=q,\ p_{n}=p$ and the set of connecting orbits $C\left(M_{\omega(s_i)}(p_{i-1}),M_{\omega(s_i)}(p_i)\right)$ is non-empty, where $\omega:[0,1]\rightarrow \Lambda$ is a path that continues $M_\lambda$ to $M_\mu$.
\end{theorem}
\prooff
By Corollary \ref{indep} a singular transition matrix is related to a generalized topological transition matrix by 
$$T_{\lambda,\mu}(p,q)=T_s\circ\times\circ\otimes\sigma(p,q)$$
therefore $T_{\lambda,\mu}(p,q)\neq 0$ imply $T_s\neq 0$. Since $T_s$ comes from a connection matrix $\Delta_g$ for a flow-defined order it follows that a non zero entry in $T_s$ implies the result. \cqd

\section{Directional Transition Matrix}

In this section,  we prove that Direction Transition Matrix is a generalized transition matrix which covers an isomorphism defined via the flow-defined Conley index isomorphism $F$. 

Consider the fast-slow systems of the form

\begin{equation}\label{fsg}
\begin{array}{ccc}
\dot{x} &=& f(x,y)\\
\dot{y} &=& \epsilon g(x,y),
\end{array}
\end{equation}

where $x\in\mathbb{R}^2$ and $y\in \mathbb{R}$. Observe that this fast-slow systems is more general than the other used to define singular transition matrix. In contrast to equation \ref{eq_fss} the slow variable $y$ depends also on the fast variable $x$, thus for $\epsilon>0$ Morse sets may have different directions of slow flow. Assume that when $\epsilon=0$ the parameterized system has an isolated invariant set $S_y$ continues over $[0,1]$ for each $y\in[0,1]$ and its Morse decomposition
$$
\mathcal{M}_y=\{M_y(\pi)\ |\ \pi\in\PP\}
$$ also continues over $[0,1]$. Furthermore, suppose that $g(M_y(\pi),y)\neq 0$ for all $y\in(0,1)$ and $\pi\in\PP$. Note that depending on $g$ the slow dynamics introduced when $\epsilon>0$ may not go in the same direction for each Morse sets. For this setting, we define an isolating neighborhood as follows.

\begin{definition}
A set $\mathcal{B}$ is a box if:
\begin{enumerate}
\item There exists an isolated neighborhood $\mathcal{B}\subseteq \mathbb{R}^n\times [0,1]$ for the parameterized flow $\psi^\mathcal{B}$ defined by
$$
\begin{array}{rcl}
\psi^\mathcal{B}:\mathbb{R}\times\mathbb{R}^n\times [0,1] &\rightarrow & \mathbb{R}^n \times [0,1]\\
(t,x,y) &\mapsto & (\psi_y(t,x),y),
\end{array}
$$
where $\psi_y$ is the flow of $\dot{x}=f(x,y)$ with fixed $y$.

\item Let $S(\mathcal{B}=Inv(\mathcal{B},\psi^\mathcal{B})$. There is a Morse decomposition 
$$
\mathcal{M}(S(\mathcal{B}))=\{M(p,\mathcal{B})\ |\ p=1, \ldots, P_\mathcal{B}\},
$$
with the usual ordering on the integers as the admissible ordering. Let $\mathcal{B}_y=\mathcal{B}\cap(\mathbb{R}^n\times \{y\}),\ S_y(\mathcal{B})=Inv(\mathcal{B}_y,\psi_y)$ and let $\{M_y(p,\mathcal{B}\ |\ p=1, \ldots, P_\mathcal{B}\}$ be the corresponding Morse decomposition of $S_y(\mathcal{B})$. Then
$$
S_0(\mathcal{B})=\displaystyle\bigcup_{p=1}^{P_\mathcal{B}}M_0(p,\mathcal{B})\ \ \ \text{and}\ \ \ S_1(\mathcal{B})=\displaystyle\bigcup_{p=1}^{P_\mathcal{B}}M_1(p,\mathcal{B}).
$$

\item There are isolating neighborhood $V(p,\mathcal{B})$ for $M(p,\mathcal{B})$ such that 
$$
V(p,\mathcal{B})\subseteq \mathcal{B}\ \ \ \text{and} \ \ \ V(p,\mathcal{B})\cap V(q,\mathcal{B})=\emptyset
$$
for $p\neq q$ with $p,q=1, \ldots, P_\mathcal{B}$, and for every $y\in [0,1]$
$$
V_y(p,\mathcal{B}\subseteq Int(\mathcal{B}_y).
$$
Furthermore, there are $\delta(p,\mathcal{B})\in\{-1,1\}, p=1, \ldots, P_\mathcal{B}),$ such that 
$$
\delta(p,\mathcal{B})g(x,y)>0,\ \ \ \text{for all}\  (x,y)\in V(p,\mathcal{B}).
$$
\end{enumerate}
\end{definition}

From the last property, one can decompose the finite index set of the Morse decomposition as 
$$\PP=\PP_+\cup \PP_-$$
where
$$\PP_\pm=\{p\in \PP\ |\ \pm\delta(p)>0\},$$
and correspondingly, one can define $M_{in}(p,\mathcal{B})$ and $M_{out}(p,\mathcal{B})$ as follows:
$$
M_{in}(p,\mathcal{B})=
\left\{
\begin{array}{cc}
M_0(p,\mathcal{B}) & \text{if}\ p\in \PP_+,\\
M_1(p,\mathcal{B}) & \text{if}\ p\in \PP_-;
\end{array}\right.
$$
$$
M_{out}(p,\mathcal{B})=
\left\{
\begin{array}{cc}
M_1(p,\mathcal{B}) & \text{if}\ p\in \PP_+,\\
M_0(p,\mathcal{B}) & \text{if}\ p\in \PP_-;
\end{array}\right.
$$

Notice that there are no connecting orbits among the Morse sets at $y=0$ and at $y=1$ and by the construction the sets $S_0(\mathcal{B})$ and $S_0(\mathcal{B})$ are related by continuation. A box with bidirectional slow dynamics can naturally occur in various problems, for instance, in the FitzHugh-Nagumo equation.  See \cite{GKMOR} and \cite{Ku} for more explanation. For this situation, either singular or topological transition matrix is not useful since they are both essentially unidirectional.

\begin{proposition}\label{contas}
Let $V$, $V'$ and $W$, $W'$ be mutually isomorphic finitely generated free Abelian groups, and let 
$$
A:V\otimes W \rightarrow V'\otimes W'
$$
be an isomorphism. Suppose $A$ is an upper triangular with the following block decomposition 
$$
A=\left( 
\begin{array}{cc}
X&Y\\
0&Z
\end{array}\right) 
$$
where $X:V\rightarrow V'$ and $Z:W\rightarrow W'$ are isomorphisms, then the following maps are all upper triangular isomorphisms:
$$
A_1=\left( 
\begin{array}{cc}
X&YZ^{-1}\\
0&Z^{-1}
\end{array}\right) : V\otimes W' \rightarrow V'\otimes W,
$$
$$
A_2=\left( 
\begin{array}{cc}
X^{-1}&-X^{-1}Y\\
0&Z
\end{array}\right) : V'\otimes W \rightarrow V\otimes W',
$$
$$
A_3=\left( 
\begin{array}{cc}
X^{-1}&-X^{-1}YZ^{-1}\\
0&Z^{-1}
\end{array}\right) : V'\otimes W' \rightarrow V\otimes W.
$$
\end{proposition}

In order to have a map from $M_{out}$ to $M_{in}$, one can repeatedly apply the Proposition \ref{contas} to the topological transition matrix until obtain an isomorphism
$$
D:\bigoplus_{p\in\PP}CH_\ast(M_{out}(p))\rightarrow \bigoplus_{p\in\PP} CH_\ast (M_{in}(p)).
$$
The matrix representation of this isomorphism is called by \textit{directional transition matrix}, which has the following property.

\begin{theorem}\emph{[KMO]}
Let $D$ be the directional transition matrix for a box in the fast-slow system (\ref{fsg}). If its $(p,q)$-entry $D(p,q)$ is nonzero, then there exist a finite sequence $\{y_i\}_{i=1}^{k+1}$ in $[0,1]$ and a sequence $\{p_i\}$ in $\PP$ satisfying
$$
\partial(p_{i+1})(y_{i+1}-y_i)>0\ \text{for all}\ i=1, \ldots, k-1
$$
and
$$
p=p_1>p_2>\ldots>p_k>p_{k+1}=q
$$
such that the corresponding parameterized flow at $y=y_i$ has a connecting orbit from $M_{y_i}(p_i)$ to $M_{y_i}(p_{i+1}).$
\end{theorem}

Simple examples show us that $D$ depends on the choices of applying Proposition \ref{contas} repeatedly and therefore led us to the question: which way should we obtain $D$ from $T$? Theorem \ref{Dcovers} answers this question by propounding that different $D$ (obtained from same $T$) cover different isomorphim, which is releated by the choice made by applying Proposition \ref{contas} repeatedly.

In this sense, the definition of directional transition matrix $D$ is deeper than just a rearrangement Morse sets from a topological transition matrix. At first glance, one could think that is an artificial definition, however $D$ is actually a transition matrix which covers an isomorphism. 

The next proposition seems to be a redundant way to define $D$, nevertheless this new way is really helpful in the continuation context as one can see in Theorem \ref{Dcovers}

\begin{proposition}\label{otherD}
A directional transition matrix $D:\bigoplus_{p\in\PP}CH_\ast(M_{out}(p))\rightarrow \bigoplus_{p\in\PP} CH_\ast (M_{in}(p))$ can be represented by an isomorphism $\bar{D}:\bigoplus_{p\in\PP}CH_\ast(M_{1}(p))\rightarrow \bigoplus_{p\in\PP} CH_\ast (M_{0}(p))$ after doing a changing of base, moreover $D$ and $\bar{D}$ are represented by the same matrix.
\end{proposition}
\prooff
Recall $$
M_{in}(p,\mathcal{B})=
\left\{
\begin{array}{cc}
M_0(p,\mathcal{B}) & \text{if}\ p\in \PP_+,\\
M_1(p,\mathcal{B}) & \text{if}\ p\in \PP_-,
\end{array}\right.
\ \ \text{and}\ \ \ 
M_{out}(p,\mathcal{B})=
\left\{
\begin{array}{cc}
M_1(p,\mathcal{B}) & \text{if}\ p\in \PP_+,\\
M_0(p,\mathcal{B}) & \text{if}\ p\in \PP_-,
\end{array}\right.
$$ define 
$$
R_{1,out}=
\left\{
\begin{array}{cc}
id(p) & \text{if}\ M_{out}(p)=M_1(p),\\
F_{01}(p) & \text{if}\ M_{out}(p)=M_0(p),
\end{array}\right.
\ \ \text{and}\ \ \ 
R_{in,0}=
\left\{
\begin{array}{cc}
id(p) & \text{if}\ M_{in}(p)=M_0(p),\\
F_{10}(p) & \text{if}\ M_{in}(p)=M_1(p).
\end{array}\right.
$$

Note that $D$ is defined via a topological transition matrix, which is an isomorphism from a base $\mathcal{B}_0$ to the base $\bigoplus_{p\in\PP}F_{01}(p)(\mathcal{B}_0)$. Thus $\bar{D}=R_{in,0}\circ D\circ R_{1,out}: \bigoplus_{p\in\PP}CH_\ast(M_{1}(p))\rightarrow \bigoplus_{p\in\PP} CH_\ast (M_{0}(p))$ is just a change of base since $R_{1,out}^{-1}=R_{in,0}$. Observe that $id$ covers $\bigoplus_{p\in\PP}F_{01}(p)$, in other words,  $\bigoplus_{p\in\PP}F_{01}(p)$ does not give information about connecting orbits between Morse sets, therefore $\bar{D}$ and $D$ are represented by the same matrix. \cqd

Same idea works, in Proposition \ref{otherD}, when one needs to change the map $\bar{D}:\bigoplus_{p\in\PP}CH_\ast(M_{1}(p))\rightarrow \bigoplus_{p\in\PP} CH_\ast (M_{0}(p))$ for $\bar{D}':\bigoplus_{p\in\PP}CH_\ast(M_{0}(p))\rightarrow \bigoplus_{p\in\PP} CH_\ast (M_{1}(p))$ by choosing another path orientation in the parameter space $[0,1]$. This can happen in Theorem \ref{Dcovers} when $n$ is even.

Applying Proposition \ref{otherD} on Proposition \ref{contas} in the fast-slow systems setting, the next lemma gives us what kind of isomorphisms that the matrices in Proposition \ref{contas} covers.

\begin{lemma}\label{contascover}
Let 
$$
T=\left( 
\begin{array}{cc}
X&Y\\
0&Z
\end{array}\right)
$$
be the topological transition matrix for the fast-slow systems in (\ref{fsg}) when $\epsilon=0$. For $A=T$ in Proposition \ref{contas} the matrices $A_1$, $A_2$ and $A_3$ are transition matrices which cover $$\left(\bigoplus_{p\in\I} F_{10}(p)\oplus F_{10}(\J)\right)\circ F_{01}(\PP)\circ \left(\bigoplus_{p\in\I} F_{10}(p)\oplus F_{10}(\J)\right),$$
$$\left(F_{10}(\I) \circ\bigoplus_{p\in\J} F_{10}(p)\right)\circ F_{01}(\PP)\circ \left(F_{10}(\I) \circ\bigoplus_{p\in\J} F_{10}(p)\right),$$
and $F_{10}$ respectively, where $\I$ and $\J$ are interval such that $\PP=\I\J$, $X:\bigoplus_{p\in\I}CH_\ast(M_{0}(p))\rightarrow \bigoplus_{p\in\I} CH_\ast (M_{1}(p))$ and $Z:\bigoplus_{p\in\J}CH_\ast(M_{0}(p))\rightarrow \bigoplus_{p\in\J} CH_\ast (M_{1}(p))$.
\end{lemma}
\prooff
Apply the same change of base done in Proposition \ref{otherD} to the matrices $A_1$, $A_2$ and $A_3$, furthermore note that the connection matrices on parameter $0$ and $1$ are equal to zero. Thus $A_1$, $A_2$ and $A_3$ are chain maps.

Since $A_3=T^{-1}$ so $A_3$ is a transition matrix which covers $F_{10}$. For the others, firstly, let $\alpha: [0,1]\rightarrow [0,1]$ be a path such that $\alpha_{|_{[0,1/3]}}$ and $\alpha_{|_{[2/3,1]}}$ are subpath from the parameter $1$ to $0$ and $\alpha_{|_{[1/3,2/3]}}$ is a subpath from $0$ to $1$.
Observe that 
$$
A_1=\left( 
\begin{array}{cc}
X&YZ^{-1}\\
0&Z^{-1}
\end{array}\right) =\left( \begin{array}{cc}
id&0\\
0&Z^{-1}
\end{array}\right) \left( \begin{array}{cc}
X&Y\\
0&Z
\end{array}\right) \left( \begin{array}{cc}
id&0\\
0&Z^{-1}
\end{array}\right),
$$
$$
A_2=\left( 
\begin{array}{cc}
X^{-1}&-X^{-1}Y\\
0&Z
\end{array}\right)=
\left( 
\begin{array}{cc}
X^{-1}&-X^{-1}Y\\
0&Z
\end{array}\right) 
\left( 
\begin{array}{cc}
X^{-1}&-X^{-1}Y\\
0&Z
\end{array}\right) 
\left( 
\begin{array}{cc}
X^{-1}&-X^{-1}Y\\
0&Z
\end{array}\right).
$$
Thus, for the path $\alpha$, $A_1$ and $A_2$ cover 
$$\left(\bigoplus_{p\in\I} F_{10}(p)\oplus F_{10}(\J)\right)\circ F_{01}(\PP)\circ \left(\bigoplus_{p\in\I} F_{10}(p)\oplus F_{10}(\J)\right),$$
$$\left(F_{10}(\I) \circ\bigoplus_{p\in\J} F_{10}(p)\right)\circ F_{01}(\PP)\circ \left(F_{10}(\I) \circ\bigoplus_{p\in\J} F_{10}(p)\right),$$
respectively, since $X$, $X^{-1}$, $Z$, $Z^{-1}$, $id(\I)$ and $id(\J)$ are transition matrices which cover $F_{01}(\I)$, $F_{10}(\I)$, $F_{01}(\J)$, $F_{10}(\J)$, $\bigoplus_{p\in\I} F_{10}(p)$ and $\bigoplus_{p\in\J} F_{10}(p)$, respectively. \cqd

The next theorem describes how direction transition matrix fits in transition matrix theory. Furthermore it shows us the relation between the choices made by applying Proposition \ref{contas} on $T$ to obtain $D$ and which isomorphism $D$ must cover.

\begin{theorem}\label{Dcovers}
Directional transition matrix is a generalized transition matrix which covers
$$
G_n\circ\cdots\circ G_1\circ F_{01} \circ G_1 \circ \cdots \circ G_n,
$$
where $G_i$ is the isomorphism defined in Lemma \ref{contascover} after applying it $i$ times, and $n$ is the number of time needed to apply Proposition \ref{contas} on topological transition matrix $T$ in order to obtain direction transition matrix $D$.
\end{theorem}
\prooff
Let $T$ be the topological transition matrix for a box $\mathcal{B}$ for a fast-slow system \ref{fsg}. Suppose that one needs to apply $n$ times Proposition \ref{contas} on $T$ in order to obtain directional transition matrix $D$. Instead of using Proposition \ref{contas}, one can use Lemma \ref{contascover}, but be aware that what was done by using Proposition \ref{contas} must be done in the same way for Lemma \ref{contascover}. In $i$-th time that one applies Lemma \ref{contascover}, the new matrix will covers $$
G_i\circ\cdots\circ G_1\circ F_{01} \circ G_1 \circ \cdots \circ G_i,
$$ thus the process ends for $i=n$. And to recover $D$, one just needs to apply  Proposition \ref{otherD}. \cqd

\end{document}